\documentclass[a4paper,11pt]{article}
\usepackage{geometry,amssymb,amsmath,latexsym}
\usepackage{xypic,xspace,enumerate}
\xyoption{curve}
\newcommand{\A}{\alpha}
\newcommand{\ba}{\bar{\alpha} }
\def\bsi{\bar{\sigma}}
\newcommand{\Ker}{\mbox{Ker}\; }

\newcommand{\sast}{_{\ast}}
\newcommand{\Env}[2]{\begin{#1} #2 \end{#1}}
\newcommand{\Sets}{\mathsf{Set}}
\newcommand{\Groups}{\mathsf{Group}}
\newcommand{\DirG}{\mathsf{Graph}}
\newcommand{\Gpds}{\mathsf{Gpd}}
\newcommand{\XMod}{\mathsf{XMod}}
\newcommand{\PXMod}{\mathsf{PXMod}}
\newcommand{\Crs}{\mathsf{Crs}}
\newcommand{\Cov}{\mathsf{Cov}}
\newcommand{\eps}{\varepsilon}
\newcommand{\io}{^{-1}}

\def\le{\leqslant}
\def\ge{\geqslant}
\def\geq{\geqslant}
\def\ow{\omega}

\newcommand{\llabto}[2]{\stackrel{#2}
{\rule[0.5ex]{#1 em}{0.07ex}\hspace{-0.4em}\longrightarrow}}

\newcommand{\mto}{\mapsto}

\def\phi{\varphi}
\newcommand{\wt}{\widetilde}
\def\wtX{\wt{X}}
\def\wtR{\wt{R}}
\def\wtG{\wt{G}}
\def\wtd{\wt{\delta}}
\def\wtI{\wt{I}}
\def\wtJ{\wt{J}}

\def\wtC{\wt{C}}
\def\wtN{\wt{N}}
\def\wtF{\wt{F}}
\def\kk{h_2'}

\newtheorem{example}{Example}[section]
\newtheorem{Def}[example]{Definition}
\newtheorem{rem}[example]{Remark}
\newtheorem{prop}[example]{Proposition}
\newtheorem{cor}[example]{Corollary}
\newtheorem{thm}[example]{Theorem}

\newtheorem{blank}[example]{\hspace{-0.45em}}

\newenvironment{proof}{\noindent {\bf Proof} }{ \hfill
$\Box$ \mbox{}}

\newcommand{\cP}{{\mathcal P}}

\def\Z{\mathbb{Z}}
\def\cc{ \mathsf{C}}

\newcommand{\sqdiagram}[8]{ \diagram  #1  \rto^-{#2} \dto_{#4} &
#3  \dto^{#5} \\ #6    \rto_-{#7}  &  #8   \enddiagram }
\def\lan{\langle}
\def\ran{\rangle}
\def\pt{\partial}
\def\mto{\mapsto}

\begin{document}
\title{\Large \bf Free crossed resolutions of groups and \\
                    presentations of modules of identities among relations}
\author{ RONALD BROWN \thanks{Part of this work was supported by EPSRC Grant
GR/J63552 `Non abelian homological algebra',1994-6, awarded to R. Brown
and T. Porter.  }\\ {\footnotesize \it School of Mathematics,
University of Wales,
 Dean Street, Bangor, Gwynedd LL57 1UT, U.K.}\\
 ABDUL RAZAK  SALLEH\thanks{This author thanks the Association
of Commonwealth Universities for support as a Commonwealth Fellow in
1994-95,  the
School of Mathematics, University of Wales, Bangor,
 for hospitality in this period,  and Universiti Kebangsaan
  Malaysia, for granting the sabbatical leave.} \\
{\footnotesize \it Department of Mathematics, Faculty of Mathematical
 Sciences,} \\
{\footnotesize \it Universiti Kebangsaan Malaysia, 43600 Bangi,
 Selangor D.E., MALAYSIA}}
\maketitle

\begin{abstract}
We give formulae for a module presentation of the module of
identities among relations for a presentation of a group, in terms
of information on 0- and 1-combings of the Cayley graph. This is
seen as a special case of extending a partial free crossed resolution of
a group given a partial contracting homotopy of its universal
cover.
\end{abstract}

\vspace{2in}


\noindent {\sc AMS Classification}:
\section*{Introduction}

The initial motivation for this work was to determine
algebraically a presentation for the $G$-module $\pi(\cP)$ of
identities among relations for a presentation $\cP= \lan X|\ow:R
\to F(X)\ran$ of a group $G$. Here we regard $R$ as a set disjoint
from $F(X)$ and $\ow$ gives the corresponding element of $F(X)$.
Recall that $\pi(\cP)$ is given algebraically as the kernel of
$\delta_2: C(R) \to F(X)$, the free crossed module of the
presentation, and is given geometrically as $\pi_2(K(\cP))$, the
second homotopy group of the cell complex of the presentation.

Our main results imply a formula as follows:

\noindent {\bf Theorem A} {\em The module $\pi(\cP)$ is generated
as $G$-module by elements
$$ \delta_3[g,r]=   (k_1 (g,\ow r))\io  \;
r ^{(\sigma g)\io }$$
for all $g \in G, r \in R$, where (i) $ \sigma: G \to F(X)$ is a
section of the quotient mapping $\phi : F(X) \to G$, (ii) $k_1$
is a morphism $F(\wt{X}) \to C(R)$ from the free groupoid on
$\wt{X}$,  the Cayley graph  of the presentation, to the free crossed
 module of the presentation, such that
  $\delta_2k_1(g,x)= (\sigma g)  x (\sigma (g(\phi x)))\io$,
  for all $x \in X, g \in G$.  }

The identities $\delta_3[g,r]$ may be seen as separation elements in
the geometry of the {\em Cayley graph with relators}, as defined
in sections \ref{compute},\ref{presgpds}. The main feature of the
theorem is that these elements generate all identities, since it
is easy to see from properties (i), (ii) and the first crossed module
rule  that these elements are all identities.

The identities $\delta_3[g,r]$ will be seen to arise from a
boundary mapping $\delta_3: C_3(I) \to C(R)$ from the free
$G$-module on a set $I$ bijective with $G \times R$, with basis
elements written $[g,r], g \in G, r \in R$. The set $\delta_3(I)$
is  usually not a minimal set of generators (many of them may even
be trivial). So we suppose given  a subset $J$ of $I$, determining
a free $G$-module $C_3(J)$, and minimal with respect to the
property that $\delta_3(J)$ also generates $\pi(\cP)$, and then
seek relations among these generators $\delta_3(J)$.

\noindent {\bf Theorem B} {\em A $G$-module generating set of
relations among these generators  $\delta_3(J)$ of $\pi(\cP)$ is
given by $$ \delta_4[g,\gamma]= -k_{2} (g,\delta_3 \gamma)
+\gamma. g\io $$ for all $g \in G, \gamma \in J$, where
$k_2: C(\wt{R}) \to C_3(J) $ is a morphism from the free crossed
$F(\wt{X})$-module on $\wt{\delta}_2: G \times R \to F(\wt{X})$
such that $k_2$ kills the operation of $F(\wtX)$ and is determined
by a choice of writing the  generators $\delta_3[g,r] \in
\delta_3(I)$   for $\pi(\cP)$ in terms of the elements of
$\delta_3 (J)$. }

It will be noted that both these results use the language of
groupoids which is convenient for encoding the graphical information.
We use in an essential way morphisms from a groupoid to a group.

 In
section \ref{compute} we shall explain the terms in these theorems
in sufficient detail for the reader to follow an explicit
calculation for the standard presentation of the group $S_3$ in
section \ref{syzygs3}. We give this example because it is
sufficiently complex to illustrate important features of the
calculations, and sufficiently simple that the calculations can be
carried out by hand.

In this example, Theorem A gives 18 generators for the module
$\pi(\cP)$;  we show this number can be reduced to 4. This
minimal set of generators   was already known. The rewriting
involved in this reduction process is then used to construct the
next level of syzygies, using Theorem B. This yields initially 24
relations among identities which are then shown to reduce to 5
independent ones. We are not aware of any previous determination
of the relations among these identities. These calculations have
been extended by hand, but with different choices, by two further
stages in \cite{He}.

The reader will notice the analogy between the formulae in these
theorems -- they are in fact special cases of Corollary
\ref{syzygies}, which computes higher order syzygies inductively.
The context of that result is that of free crossed resolutions,
universal covers of crossed resolutions, and contracting
homotopies of such universal covers. Once this machinery is set
up, the result becomes almost tautologous. It states that the pair
consisting of a partial free crossed resolution and a partial
contracting homotopy of its universal cover can be extended by one
step, and hence indefinitely.

A sequel to this paper by Heyworth and Wensley \cite{HeW} will
show how the  part of the procedure required for Theorem A can be
implemented as a `logged Knuth-Bendix procedure'. A further
paper by Heyworth and Reinert \cite{HeR} will show how generalised
Gr\"{o}bner basis procedures for integral group rings can
implement the reduction process required for Theorem B,  and so
allow a wide range of computations.

The partial contracting homotopies are given by functions $h_i$
for $i < n$ with appropriate properties. In fact $h_0$ corresponds
to a 0-combing, and $h_1$ is analogous to a 1-combing; from these
we obtain the functions $\sigma, k_1$ of the theorems. The algebra
of such functions is shown to be nicely handled in the context of
the free groupoid $F(\wt{X})$ on the Cayley graph and the free
crossed module $\wt{\delta}_2: C(\wt{R}) \to F(\wt{X})$. We show
that this crossed module is the fundamental crossed module of the
universal cover of the geometric 2-complex of the presentation.
The groupoid approach is required to utilise all the vertices of
the universal cover.

A computation of the module of identities among relations for the
presentation $\lan X|R\ran $ of the group $G$ could be seen in the
context of chain complexes and resolutions of modules  as that of
computing an extension of the partial resolution of $\Z$ $$ (\Z
G)^R \to  (\Z G)^X \to \Z G \to  \Z$$ where the first morphism is
given by the Whitehead-Fox derivative $(\partial r/ \partial x)$
\cite[Lemma 8]{W1}, \cite{Fox}. The process of extending
a partial resolution is more difficult than that of just giving a
resolution. There is in fact considerable work on constructing
resolutions of groups, some of it for 2-groups mod 2, and other
results using homological perturbation theory, particularly by
Larry Lambe and colleagues \cite{JLS}. It is not clear how  these
methods apply to the problem of extending partial
resolutions. Work of Groves \cite{Gro} constructs a resolution
from a complete rewrite system for a monoid presentation of the
group, rather than directly from a group presentation. However, as
mentioned above, complete rewrite systems are relevant to the
computation of $k_1$.

It is interesting to compare our methods with the methods of
pictures for calculating the generators of $\pi(\cP)$ (see for
example \cite{B-Hu,HMS,P}). These methods use nicely the geometry
of the relations, they  have been very successful in this field,
and can be  more efficient than ours for this dimension. However
they seem more difficult to carry out in higher dimensions, for
the following reasons.

The picture methods use 2-dimensional rewrite rules to reduce
spherical elements to a combination of standard elements. The full
information on the way these rewrites are used in a particular
example is essentially 3-dimensional, and it  can thus be
difficult to visualise or to record combinations of such rewrites,
and their dependencies. For our purposes this   rewriting
information must be recorded completely (see Tables 2,4 of section
\ref{syzygs3})  since it is  used to construct the next stage of a
contracting homotopy; this use of the complete record is one
reason for the apparently cumbersome nature of the calculations.
Thus there are problems in extending the  picture method to
determine 3-syzygies, whereas our purely algebraic method is
essentially uniform over dimensions, giving rise mainly to
computational problems. This suggests that in dimension 2 our
methods  should be seen as complementary to those of pictures.

The method of pictures has also been applied successfully to
determine generators for the module of identities among relations
for various constructions on groups. By contrast the only general
construction on crossed resolutions which has so far been applied
is the tensor product \cite{Br-Hi2,B-P,T} -- given  free crossed
resolutions $C,C'$ of two groups $G,H$, the tensor product $C
\otimes C'$ gives a free crossed resolution of their  product $G
\times H$, and so a presentation of the module of identities for
the standard combined presentation of the product. An application
is in \cite{B-P}.

There are three basic planks  in our approach.

\noindent (i) {\em Crossed complexes}
\nopagebreak

Crossed complexes  form an analogue of chain  complexes but with
non abelian features in dimensions 1 and 2. These features allow
crossed complexes to combine many of the advantages of chain
complexes with an ability to contain the information involved in a
presentation of a group. So one can model many of the standard
techniques of homological algebra, such as uniqueness up to
homotopy of a free crossed resolution. Further, this technique may
be combined with a non abelian version of the traditional notion
of `chains of syzygies'; this version takes account of the facts
that free groups are non abelian, and that a normal subgroup $N$
of a group $F$  is  in general non abelian, and admits an
operation of $F$ on $N$ which is crucial in discussing
presentations. Crossed complexes, unlike chain complexes, allow
for `free' models of this inclusion $N \to F$ (see \cite{B-Hu}),
and so give an intuitive algebraic model of chains of syzygies in
this non abelian case.   An account of uses of
crossed complexes up to 1981   is given in \cite{BH-81}.

A small free crossed resolution is convenient for calculations of
non abelian extensions \cite{B-P} and of the cohomology class of a
crossed module \cite{BW,BW2}. A free crossed resolution $C$ of $G$
determines a free $\mathbb{Z}G$ resolution  $\Delta C$ of
$\mathbb{Z}$ in the usual sense \cite{W1,BH-91}. The crossed
resolution $C$ with its free basis carries more information than
$\Delta C$, for example it includes a presentation of $G$.

\noindent (ii) {\em Algebraic models of the geometry of
covering spaces}

Philip Higgins pointed out in 1964 \cite{Higg-Pres}  how
presentations of  groupoids could be applied to group theory. The
geometric basis of the argument is that the theory  of covering
spaces is more conveniently handled if one uses groupoids rather
than groups, since there is a purely algebraic notion of {\em
covering morphism of a groupoid} which nicely models the geometry
(see \cite{B-88}). Covering morphisms of a group or groupoid $G$
are equivalent to operations of $G$ on sets.

In the same way, to apply crossed complexes  to covering spaces we
require crossed complexes of groupoids not just of groups. Such
general crossed complexes were also found essential in  \cite{BH3}
for certain higher order Van Kampen Theorems, so the basic
definitions and applications are already known. This allows us to
bring in techniques not only  of presentations of groupoids, as
discussed in \cite{Higg-Pres}, but also of free crossed
resolutions determined by such a presentation.

In effect, we are giving a suitable {\em algebraic} framework in
which to place the geometry of the Cayley graph of a generating
set of the group, but including the relations as well as the
generators of the presentation, and indeed including higher order
syzygies, as these are constructed.  This algebraic framework also
models conveniently the geometry of the universal cover of a cell
complex.

A crucial tool for our methods is the fact that a covering crossed
complex of a free crossed complex is again a free crossed complex,
on the `covering generators' (Theorem \ref{T:freecov}). This
models the geometry of CW-complexes. The result is crucial because
it enables us to define morphisms and homotopies  by their values
on the free generators. Our proof relies on a result of Howie
\cite{Howie}.

\noindent (iii) {\em Contracting homotopies}

The key point is that the previous techniques  allow us to discuss
free crossed resolutions of contractible groupoids, for example
the universal covering groupoid of the original group.  A crossed
resolution of a contractible groupoid  will have contracting
homotopies, and our method proceeds by the construction of such
homotopies. This method is applied to truncated  crossed complexes
and in particular to crossed modules. The usual slogan {\em choose
generators for the kernel} and so kill homotopy groups, fails to
tell us how to choose these generators. Instead  we construct a
crossed  complex whose universal cover is a {\em home for a
contracting homotopy}. This `tautologously'  yields generators of
kernels.

In order to make this method clear, we need the basics of the
theory of presentations and of identities among relations for
groupoids. We give the key features, largely without proofs, in
section \ref{presgpds}.

The basic theory of crossed complexes and their covering morphisms
that we need is presented in sections \ref{cov}--\ref{free}.
Finally, the notion of homotopy for crossed complexes is presented
in section \ref{homotopies}.

Our method  yields a resolution dependent functorially on the
presentation. However a count of the numbers of generators in
various dimensions shows that the module resolution obtained from
our crossed resolution by the process of  \cite{BH-90} is not the
same as the Gruenberg resolution \cite{Gru1}. We are grateful to
Justin Smith for pointing out this reference.

More generally, we can obtain a free crossed resolution dependent
functorially on the first $n$ stages of a free crossed resolution,
with basis up to this stage.

In the final section we show how these methods give rise to the
standard crossed resolution of a group $G$, and to a small crossed
resolution of a finite cyclic group. In each case, the information
on the contracting homotopy determines the resolution.

We would like to thank Anne Heyworth and Emma Moore for
discussions on this material which led to the exposition in
sections \ref{compute}, \ref{moreres}.

\section{The computational procedure}
\label{compute}

The purpose of this section is to state the computational
procedure in as direct a way as we can. The theoretical underpinning
is left to later sections. We hope this will make it easier for the
reader.

Let $\cP= \lan X|\ow: R \to F(X) \ran $ be a presentation of a group
$G$. The advantages of using the function $\ow$ are (i) to allow for
the possibility of repeated relations, and (ii) to distinguish
between an element  $r \in R$ and the corresponding element
$w(r)\in F(X)$. We shall be concerned with the following diagram,
in which $p_0 \wt{\phi}=\phi p_1, p_1\wtd_2=\delta_2 p_2,
p_2\wtd_3=\delta_3 p_3$. The parts of this diagram will be developed below:

\begin{equation}
\begin{split} \label{maindiag} \def\labelstyle{\textstyle}
\xymatrixcolsep{3pc}
\xymatrix {G&G&G& G \ar @/_0.75pc/[dl] _-{h_0}\\
 C_3(\wt{I})\ar [u]^-{\beta} \ar [d] _{p_3} \ar [r] ^-{\wt{\delta}_3} &
C(\wt{R})\ar [u]^-{\beta} \ar [d] _{p_2} \ar [r] ^-{\wt{\delta}_2}\ar @/^0.75pc/ [l]^{h_2} &
 F(\wt{X})\ar [u]^-{\beta} \ar [d] _{p_1} \ar [r]^-{\wt{\phi}} \ar @/^0.75pc/ [l]^{h_1}
 &\ar@<1ex> [u] ^-{\delta^0} \ar@<-1ex> [u] _-{\delta^1}\wt{G} \ar [d] _{p_0} \\
C_3(I)\ar [r] _-{\delta_3} & C(R) \ar [r] _-{\delta_2} & F(X)
\ar [r]_-{\phi}  & G}
\end{split}
\end{equation}

\begin{blank}{\em  $\phi: F(X) \to G$ is the canonical morphism from the free
group on $X$ to $G$ given by the set of generators.}
\end{blank}
\begin{blank}
{\em $\delta_2: C(R) \to F(X)$ is the free crossed
$F(X)$-module on the function  $\ow$.
}
\end{blank}
Thus the elements of $C(R)$ are `formal consequences'
$$c=\prod_{i=1}^n (r_i^{\eps_i})^{u_i}$$
where $n \ge 0,r_i \in R, \epsilon_i=\pm 1, u_i \in F(X)$,
$\delta_2(r^{\epsilon})^u=u\io (\ow r)^{\epsilon} u$, subject to
the crossed module rule $ab=ba^{\delta_2b}, a,b \in C(R)$.
For information on crossed modules, and particularly free crossed modules,
see for example \cite{B-Hu,HMS,B-gre}.

Let $N= \Ker \phi$. Then $\delta_2(C(R))=N$. Of course it is the
kernel $\pi(\cP)$ of $\delta_2$, the $G$-module of identities
among relations, that we wish to calculate. For this we require
algebraic analogues  of methods of
covering spaces, and so use the language of groupoids. Our
convention is that the product of elements (arrows) $a:g \to g', a':g' \to
g''$ in a groupoid $\Gamma$ is written $aa': g \to g''$, and $\Gamma(a)$ denotes
the object group of $\Gamma$ at $a$, i.e. the set of arrows $a \to
a$ with the induced group structure.

\begin{blank}
{\em  $p_0: \wt{G} \to G$ is the universal covering groupoid of
the group $G$. The objects of $\wt{G}$ are the elements of $G$,
and an arrow of $\wt{G}$ is a pair $(g,g')\in G \times G$ with
source $g= \delta^0(g,g')$ and target $gg'=\delta^1(g,g')$. The
projection morphism $p_0$ is given by $(g,g') \mapsto g'$. }
\end{blank}
\begin{blank}{\em
$\wt{X}$ is the { \em Cayley graph} of the pair $(G,X)$. Its objects are
the elements of $G$ and its arrows are pairs $(g,x) \in G \times
X$ with source $g=\delta^0(g,x)$ and target
$g(\phi x)=\delta^1(g,x)$, also written $\beta(g,x)$.}
\end{blank}
\begin{blank}{\em
$F(\wt{X})$ is the {\em free groupoid} on $\wt{X}$. Its objects are the
elements of $G$ and its arrows are pairs $(g,u) \in G \times F(X)$
with source $g$ and target $g(\phi u)$.  We also write
$\beta (g,u)=g(\phi u)$. The multiplication is
given by $(g,u)(g(\phi u) ,v)= (g, uv)$. The morphism $\wt{\phi}$
is given by $(g,u) \mapsto (g, \phi u)$. The morphism $p_1$ is
given by $(g,u) \mapsto u$. It maps the object  group $F(\wtX)(1)$
isomorphically to $N$.  }
\end{blank}

As we shall see in section \ref{cov}, $\wtG \to G$ is the covering
morphism corresponding to the trivial subgroup of $G$, and
$F(\wtX) \to F(X)$ is the covering morphism corresponding to the
subgroup $N$ of $F(X)$.
\begin{blank} \label{frecovxm}{\em
$\wt{R}= G \times R$ and $\wt{\delta}_2: C(\wt{R})\to F(\wt{X})$ is
the {\em free crossed $F(\wt{X})$-module} on
$\wt{\ow}: \wt{R} \to F(\wt{X}), \;
(g,r) \mapsto (g, \ow(r))$. Then  $C(\wtR)$ is the disjoint union of
groups $C(\wtR)(g), g \in G$, all mapped by $p_2$ isomorphically to
$C(R)$. Elements of $C(\wtR)(g)$ are pairs $(g,c) \in
 \{g \}  \times C(R)$, with
multiplication $(g,c)(g,c')=(g,cc')$. The (partial) action of $F(\wt{X})$ is
given by $(g,c)^{(g,u)}= (g(\phi u), c^u)$. The boundary
$\wt{\delta}_2$ is given by $(g,c) \mapsto (g, \delta_2 c)$. The
morphism $p_2: C(\wt{R}) \to C(R)$ is given by $(g,c) \mapsto c$.}
\end{blank}

If $(g,c) \in C(\wtR)(g)$ we write $\beta (g,c) = g$ ; we call $\beta$
 the {\em base point map}.
 The elements of $C(\wtR)(g)$ are also all `formal consequences'
 $$(g,c)=\prod_{i=1}^n ((g_i,r_i)^{\eps_i})^{(g_i,u_i)}
 =\prod_{i=1}^n (g,(r_i^{\eps_i})^{u_i})=(g,\prod_{i=1}^n (r_i^{\eps_i})^{u_i}) $$
where $n \ge 0,\, r_i \in R, \, \epsilon_i=\pm 1,\,  u_i \in F(X),\,
 g_i\in G,\, g_i(\phi u_i)=g $,
 subject to the crossed module rule $ab=ba^{\wtd_2b}, a,b \in
 C(\wtR)$. Here the first form of the product is useful geometrically, and the
last is useful computationally.

In effect, we are giving first a presentation $\lan \wt{X} | \wt{\ow}:\wt{R}\to
F(\wt{X})\ran$ of the groupoid $\wt{G}$ \cite{Higg-Pres}, and second the
 free crossed module corresponding to this presentation.

The proof that the construction given in \ref{frecovxm} does give
the free crossed module as claimed is given in theorem
\ref{T:freecov}.

We now construct $C_3(I)$ and its cover $C_3(\wtI)$.

\begin{blank} {\em
Let $I$ be a set in one-to-one correspondence with $G \times R$
with elements written $[g,r], g \in G, r \in R$. Let $C_3(I)$ be
the free $G$-module on $I$.}
\end{blank}
\begin{blank}{\em
Let $C_3(\wtI)$ be the free $\wtG$-module on $\wtI=G \times I$.
 Then  $C_3(\wtI)$ is the disjoint union of
abelian groups $C(\wtI)(g), g \in G$, all mapped by $p_3$ isomorphically to
$C_3(I)$. Elements of $C_3(\wtI)(g)$ are pairs $(g,i) \in
 \{g \} \times C_3(I)$  with addition $(g,i) + (g,i') =(g,i +
i')$. The  (partial) action of $\wtG$ on $C_3(\wtI)$ is given by $(g,i).{(g,g')} =
(gg', i.{g'})$.}
\end{blank}
The construction of $\delta_3$ (and hence of $\wtd_3$) requires
some choices.

\begin{blank}{\em  Choose a section $\sigma: G \to F(X)$ of $\phi$ such that
$\sigma (1) =1,$ and write $\bsi(g)=\sigma(g)\io$. Then
$\sigma $ determines a function $h_0: G \to
F(\wt{X})$ by $g \mapsto (g, \bsi g )$. Thus $h_0(g)$ is a
path $ g \to 1$ in the Cayley graph $\wt{X}$.}
\end{blank}
\begin{rem}
{\em The choice of $h_0$ is often, but not always, made by
choosing a maximal tree in the graph $\wtX$ -- such a choice is
equivalent to a choice of Schreier transversal for the subgroup
$N= \Ker \phi $ of $F(X)$. A different choice of $h_0$ is used in
subsection \ref{standres} for the standard crossed resolution.
}
\end{rem}

For each arrow $(g,x)$ of $\wt{X}$ the element
$\rho (g,x)=(h_0 g) \io (g,x) h_0(g(\phi x))$ is a loop at $1$ in $F(\wtX)$
and so is in the image of $\wtd_2$.
\begin{blank} {\em
For each arrow $(g,x)$ of $\wt{X}$ choose an element $h_1(g,x) \in
C(\wtR)(1)$ such that
$$
\wtd_2 (h_1(g,x) )= \rho (g,x) .
$$
Then $h_1$ extends uniquely to a morphism $h_1: F(\wtX) \to C(\wtR)(1)$
 such that for all arrows $(g,u)$ of $F(\wtX)$
\begin{equation}\label{bdyh1}
\wtd_2 (h_1(g,u) )= (h_0 g) \io (g,u) h_0(g(\phi u)).
\end{equation}
It follows that $h_1(h_0(g))=1, g \in G$.
}
\end{blank}
\begin{rem}
{\em The choice of $h_1$ is equivalent to choosing  a representation as
a consequence of the relators $R$ for each element
 of $N$, given as a word in the elements of $X$. There is no algorithm for
 such a choice. It will be shown in \cite{HeW} how a `logged
Knuth-Bendix procedure' will give such a choice  when  the monoid
rewrite system determined by $R$ may be completed, and that this allows for
an implementation of the determination of $h_1$. }
\end{rem}

The morphism $k_1$ of Theorem A of the Introduction is simply the
composition $p_2 h_1$.

\begin{blank}{\em Define $\delta_3: C_3(I) \to C(R)$ by
\begin{equation} \label{delta3}
\delta_3[g,r]=  p_2\left((h_1  (g,w(r)))\io \right) \;
r ^{\bsi g }.
\end{equation}}
\end{blank}
It follows from equation \eqref{bdyh1} that $\delta_2 \delta_3=0$, and so the given
values $\delta_3[g,r]$ lie in the $G$-module $\pi(\cP)$. This implies that
$\delta_3$ is well defined on $C_3(I)$ by its values on the set
$I$ of free module generators.

\begin{blank}{\em
 Let $\wtd_3: C_3(\wtI) \to C(\wtR)$ be the $\wtG$-morphism
given by
$\wtd_3(g,d)= (g, \delta_3 d)$. Let  $h_2:C(\wtR) \to
C_3(\wtI)(1)$ be the groupoid morphism killing the operation of
$F(\wtX)$ (i.e. $h_2((g,c)^{(g,u)})= h_2(g,c)$ for all $(g,c) \in C(\wtR), u \in F(X)$)
and satisfying $(g,r) \mapsto (1,[g,r]), (g,r) \in G \times R$.
Then for all $g \in G, c \in C(R) $
\begin{equation} \label{bdyh2}
\wtd_3 h_2(g,c) = (h_1\wt{\delta}_2 (g,c))\io  \; (1,c^{\bsi g }).
\end{equation}
}
 \end{blank}

\begin{blank} {\bf Proof of Theorem A} {\em
Equations \eqref{bdyh1}, \eqref{bdyh2} show that $\wtd_2 \wtd_3=0$,
and so the elements $p_2(\wtd_3 h_2(g,c))$ do give identities. On
the other hand, if $c \in C(R)$ and $\delta_2 c =1$, then by
equation \eqref{bdyh2}, $(1,c)= \wtd_3h_2(1,c)$, and so $c=\delta_3(d)$
 for some $d$.
Theorem A of the Introduction is an immediate consequence,
 with $k_1=p_2h_1$.
} \hfill $\Box$
\end{blank}

However some of the elements of $\delta_3(I)$ may be trivial, and others
may depend $\Z G$-linearly on a smaller subset. That is, there
may be a proper subset $J$ of $I$ such that $\delta_3(J)$ also
generates the module $\pi(\cP)$. Then for each element $i\in I \setminus
J$ there is a formula expressing $\delta_3 i$ as a $\Z G$-linear
combination of the elements of $\delta_3( J)$. These formulae determine
a $\Z G$-retraction $\mathsf{r} : C_3(I) \to C_3(J)$ such that for all $ d
\in C_3(I), \delta_3(\mathsf{r} d) = \delta_3(d)$. So we replace $I$ in
the above diagram by $J$, replacing the boundaries by their
restrictions. Further, and this is the crucial step, we replace
$h_2$ by $h'_2= \mathsf{r}'  h_2$ where $\mathsf{r} ': C_3(\wtI)(1) \to C_3(\wt{J})(1)$
is mapped by $p_3$ to $\mathsf{r}$.

This $h'_2: C(\wtR)\to C_3(\wt{J})(1)$
is now used to continue the above construction, by defining
$C_4(\bar{J})$ to be the free $G$-module on elements written
$[g,d] \in \bar{J}= G \times J$, with
\begin{equation} \label{2syzy}
 \delta_4 [g,d]= -p_3(h'_2 (g,\delta_3 d)) + d . g \io.
\end{equation}
These boundary elements give generators for the relations among the generators
 $\delta_3(J)$ of $\pi(\cP)$.

\begin{blank} {\bf Proof of Theorem B} {\em
This is a similar argument to the proof of Theorem A, using
 equation \eqref{2syzy}, and setting $k_2=p_3h_2'$.
} \hfill $\Box$ \end{blank}

\begin{rem}{\em
In the above we have defined morphisms and homotopies by their
values on certain generators, and so it is important for this that
the structures be free. For example, $h'_2$ is defined by its
values on the elements $(g,r) \in G \times R$. So, noting that $h_2$
kills the operation of $F(\wtX)$,  we calculate for
example $h'_2(g,r^us^v)= h'_2(g(\phi u)\io,r) + h'_2(g(\phi v)\io,s)$.
 In this way the formulae reflect the choices made at different parts of
the Cayley graph in order to obtain a contraction.

The freeness of $F(\wtX)$ was proved by Higgins
in \cite{Higg-Pres}. Our proof for $C(\wtR)$ uses a result of
Howie, as we shall see later.
}\end{rem}
\begin{rem} {\em
The determination of  minimal subsets $J$ of $I$ such that
$\delta_3 J$ also generates $\pi(\cP)$ is again not
straightforward.  Some dependencies are easy to find, and others
are not. A basic result  due to Whitehead \cite{W1} is that the
abelianisation map $C(R) \to (\Z G)^R$ maps $\pi(\cP)$
isomorphically to the kernel of the Whitehead-Fox derivative $(\pt
r/\pt x):(\Z G)^R \to (\Z G)^X$. Hence we can test for dependency
among identities by passing to the free $\Z G$-module $(\Z G)^R$, and we
use this in the next section. For bigger examples, this testing can be a formidable
task by  hand. An implementation of Gr\"{o}bner basis procedures
for finding minimal subsets which still generate is described in \cite{HeR}.
}
\end{rem}

\section{Syzygies of levels 2 and 3  for the standard presentation
of $ \bf S_3$}
\label{syzygs3}

We illustrate the above method in this section with the
standard presentation of the six element group $S_3$. This is
chosen as perhaps the first interesting example which can still be
done by hand, and because it does illustrate all the above
points. While our set of generators of the module of identities  for
this presentation is known, we are not aware of previous calculations of
the relations between these generators.

The group presentation $\langle
x,y\mid x^3,y^2,xyxy\rangle$ determines  the symmetric group  $S_3$
on three symbols. Let $X= \{x,y\}$ and let $F=F(X)$ be the free group on
$X$. Let $R= \{r,s,t\}$ and let $\ow :R \to F$ be given by $$\ow r=x^3,
\ow s=y^2,\ow t=xyxy.$$  Let $\phi : F \to S_3$ be the epimorphism
determined by $x,y$, and let $N= \Ker \phi.$ So we have the free
crossed module $\delta_2: C(R) \to F$.

Now we set up the corresponding diagram \eqref{maindiag} of the
previous section. We think of each element $(g,r) \in \wtR$ as
filling a 2-cell in the Cayley graph $\wtX$.
Thus in this example, each relator, i.e. each element of $R$, is covered
 six times in the universal cover. We also see in this situation the
r\^{o}le of relations which are proper powers. The covers of
the element $r$ of $R$ separate into two classes, namely $$\{
(1,r), (\phi x,r), (\phi x^2,r) \}, \qquad \{(\phi y,r), (\phi
yx,r), (\phi yx^2,r) \}.$$ An element of one of these  classes has
boundary the same `triangle' in $F(\wt{X})$ as the other
elements, but with a different starting point. Similarly, the
relation $\ow t$ is of order 2 and so the covers of $t$ divide
into three classes each with 2 elements. A similar statement holds
for $s$.

We now have to choose $\sigma : S_3 \to F(X)$. For this, choose a maximal
tree $T$ in the directed graph $\wt{X}$. The choice
of $T$ is well known to be equivalent to the choice of a Schreier
transversal for $N$ in $F$. For this example, we choose the tree $T$ to be
given by the elements
$$(1,y), \, (1,x), \, (\phi x^2, x), \, (\phi y, x), \, (\phi (xy),
x).$$
The remaining elements of $\wt{X}$ we label as
\begin{gather}\theta _1=(\phi y,y), \; \theta _2= (\phi x, y), \;  \theta _3= (\phi xy,
y),\; \theta _4=(\phi yx,y),\\ \theta _5= (\phi x^2, y), \;  \theta _6 =(\phi yx,
x),\; \theta _7=( \phi x,x). \end{gather}
\setlength{\unitlength}{1.5cm}%
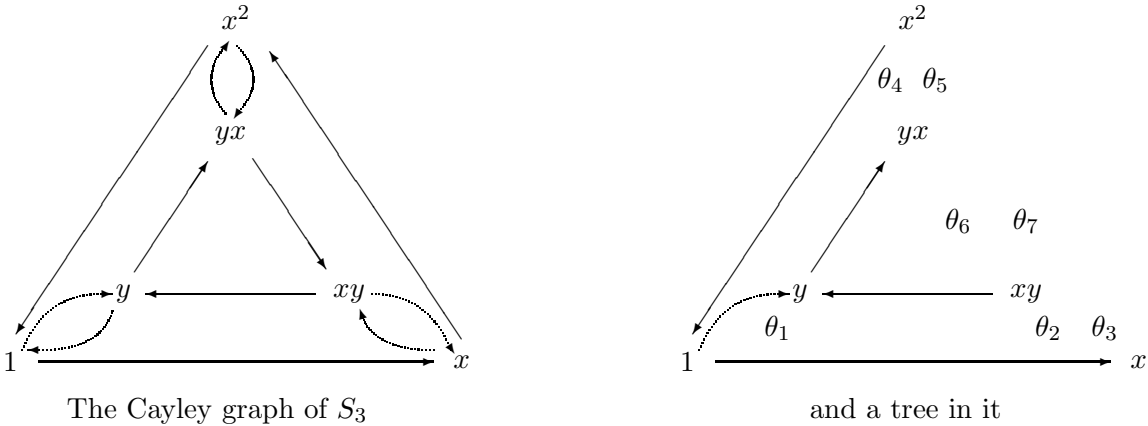
\begin{picture}(15,4.2)(-2.7,-1)
\put(-2,0){\makebox(0,0){1}}
\put(-1.75,0){\vector(1,0){3.5}}
\put(2,0){\makebox(0,0){$x$}}
\put(2,0.2){\vector(-2,3){1.7}}
\put(0,2.95){\makebox(0,0)[b]{$x^2$}}
\put(-0.25,2.8){\vector(-2,-3){1.7}}

\put(-1,0.6){\makebox(0,0){$y$}}
\put(-0.9,0.8){\vector(2,3){0.65}}
\put(0.7,0.6){\vector(-2,0){1.5}}
\put(-0.05,2){\makebox(0,0){$yx$}}
\put(0.15,1.8){\vector(2,-3){0.65}}
\put(1,0.6){\makebox(0,0){$xy$}}
\qbezier[30](-1.9,0.1)(-1.7,0.6)(-1.2,0.6)
\qbezier[30](-1.75,0.1)(-1.2,0.1)(-1.1,0.45)
\put(-1.2,0.6){\vector(1,0){0.1}}
\put(-1.75,0.1){\vector(-1,0){0.1}}
\qbezier[30](-0.1,2.2)(-0.35,2.5)(-0.1,2.8)
\put(-0.16,2.7){\vector(2,3){0.1}}
\qbezier[30](0,2.2)(0.3,2.5)(0,2.8)
\put(0.08,2.3){\vector(-2,-3){0.1}}

\qbezier[30](1.9,0.1)(1.7,0.6)(1.2,0.6)
\qbezier[30](1.75,0.1)(1.2,0.1)(1.1,0.43)
\put(1.117,0.41){\vector(-1,3){0.02}}
\put(1.84,0.215){\vector(2,-3){0.1}}

\put(4,0){\makebox(0,0){1}}
\put(4.25,0){\vector(1,0){3.5}}
\put(8,0){\makebox(0,0){$x$}}
\put(6,2.95){\makebox(0,0)[b]{$x^2$}}
\put(5.75,2.8){\vector(-2,-3){1.7}}

\put(5,0.6){\makebox(0,0){$y$}}
\put(5.1,0.8){\vector(2,3){0.65}}
\put(6.7,0.6){\vector(-2,0){1.5}}
\put(6,2){\makebox(0,0){$yx$}}

\put(7,0.6){\makebox(0,0){$xy$}}
\qbezier[30](4.1,0.1)(4.3,0.6)(4.8,0.6)
\put(4.8,0.6){\vector(1,0){0.1}}

\put(4.8,0.3){\makebox(0,0){$\theta _1$}}
\put(7.2,0.3){\makebox(0,0){$\theta _2$}}
\put(7.7,0.3){\makebox(0,0){$\theta _3$}}
\put(6.4,1.25){\makebox(0,0){$\theta _6$}}
\put(7.0,1.25){\makebox(0,0){$\theta _7$}}
\put(5.8,2.5){\makebox(0,0){$\theta _4$}}
\put(6.2,2.5){\makebox(0,0){$\theta _5$}}

\put(-1.5,-0.5){\mbox{The Cayley graph of $S_3$}}
\put(5,-0.5){\mbox{ and a tree in it}}

\end{picture}

The object groups of  the free groupoid $F(T)$ on the graph $T$
are all trivial. For each $g \in S_3$ let $h_0 g$ denote the
unique element of $F(T)(g,1)$, so that the section $\sigma
: S_3 \to F$ of $\phi$ is given by $h_0g= (g, (\sigma g) \io), g \in G$.  Then for
each $(g,u):g \to g'$ in $F(\wtX)$, we set $\rho (g,u) =
(h_0(g))\io (g, u) h_0(g') $. Let $\wtN = F(\wtX)(1)$, which is
mapped isomorphically by $p_2$ to  $N=\delta_2(C(R))=\Ker \phi$.
Thus the tree $T$ determines a retraction $\rho:F(\wtX)\to
\wt{N}$.

Let $D$ be the set of edges of $\wt{X}$ which do not lie in $T$.
Then the set  $\rho (D)$ is a set of free generators of the group
$\wt{N}$, and $p_1 \rho(D)$ is a set of free generators of the
group $N$. Let  $\eta_i = \rho \theta _i, i=1,\ldots, 7.$

In order to define $h_1: F(\wtX) \to C(\wt{R})(1)$ we need only to
give its values on the generators (see \ref{spechomot}). We give
these by $h_1(\tau)=1 $ if $ \tau \in T$ and for $\theta \in D$,
we let $h_1(\theta)$ be an element of $C(\wtR)(1)$ which is mapped
by $\wtd_2$ to $\rho(\theta)$. Then $h_1$ satisfies (\ref{bdyh1}),
and also $h_1(\theta)= h_1(\rho(\theta))$.

In our example of $S_3$, we  define $h_1$ on $\rho(D )$, and so on
$F(\wtX)$,  as follows:
\begin{alignat*}{3}
\eta_1 &=\rho (\phi y,y) & \qquad h_1 \eta _1 &= (1,s)& &=(1,s) \\ \eta_2
&= \rho (\phi x,y) &\qquad  h_1 \eta _2 &= (1,t)(1,s)\io&&= (1,ts\io) \\ \eta_3
&=\rho (\phi xy,y)  & \qquad h_1 \eta _3 &= (1,s)(1,t)\io (\phi
x,s)^{(1,x)\io}&&= (1,st\io s^{x\io}) \\ \eta_4 &= \rho (\phi yx, y)& \qquad h_1 \eta _4
&=(\phi y,t)^{(1,y) \io} &&= (1,t^{y\io})
\\ \eta_5 &= \rho (\phi x^2,y) & \qquad h_1 \eta_5 &=
(\phi x^2,s)^{(\phi x^2,x)}((\phi x^2,t)\io)^{(\phi x^2,x)} &&=(1,s^x(t\io)^x) \\ \eta_6 &= \rho (\phi yx,x)
 &\qquad  h_1 \eta
_6    &= (\phi y,r)^{(1,y)\io}&&= (1,r^{y\io}) \\ \eta_7    &=\rho  (\phi x,x) &
\qquad h_1 \eta _7    &= (1,r) &&= (1,r)
\end{alignat*}
This ensures that $\wt{\delta}_2 h_1 (\eta _i) = \eta_i, i=1,2,\ldots,7. $

In order to calculate the identities among relations, we now need to
express $\rho\wt{\delta}_2 \alpha$ in terms of the $\eta _i$ for all
$\alpha \in \wt{R}$.
Then according to the previous section we can obtain
 an  identity among relations  for each $\alpha \in \wt{R}, $
namely $$p_2\left( (h_1\wt{\delta}_2 \alpha)\io\alpha^{h_0 \beta
\alpha}\right).$$ The results of these calculations are given in
the  table which follows. The order of writing the identities is
chosen so that the first four give our
eventual minimal set of generators, the next six give trivial
identities, and the last has the most difficult  verification of
its dependence on the first four.
\begin{center}
Table 1
\nopagebreak

\begin{tabular}{||l|l||l||l|l||}
\hline \multicolumn{2}{||c||}{generator} & $\rho\wt{\delta}_2
\A_i$& \multicolumn{2}{c||}{\rule{0em}{3ex} $\gamma_i=
p_2\left(h_1(\wt{\delta}_2 \A_i)\io \A_i^{h_0 \beta
\A_i}\right)$}\\ \hline $\alpha _1$   & $(\phi x^2,r)$& $\eta_7$
&$\gamma_1$  &$r\io r^{x}$\\ $\alpha _2$   & $(\phi y, s)$ & $\eta
_1 $      &$\gamma_2$  &$s\io s^{y\io}$ \\ $\alpha _3$   & $(\phi
x^2,s)$& $\eta _5\eta_4 $&$\gamma_3$  &$(t\io)^{y\io}t^x$\\
$\alpha _4$   & $(\phi x,t)$  & $\eta_7 \eta_5\eta_6 \eta _3$&
$\gamma_4$ &
        $(s\io)^{x\io}ts\io (r\io)^{y\io}t^x(s\io)^xr\io t^{x\io}$\\
$\alpha _5$   & $(1,r)$       & $\eta_7$        &$\gamma_5$    & 1 \\
$\alpha _6$   & $(1,s)$       & $\eta _1 $      &$\gamma_6$    &1 \\
$\alpha _7$   & $(1,t)$       & $\eta_2 \eta_1$ &$\gamma_7$    &1 \\
$\alpha _8$   & $(\phi x,s)$  & $\eta_2 \eta_3 $&$\gamma_8$    & 1 \\
$\alpha _9$   & $(\phi y,t)$  & $ \eta _5 $     &$\gamma_9$    & 1 \\
$\alpha _{10}$& $(\phi y,r)$  & $ \eta_6$       &$\gamma_{10}$ &$1$\\
$\alpha _{11}$& $(\phi x,r)$  & $\eta_7 $       &$\gamma_{11}$ &$r\io r^{x\io}$ \\
$\alpha _{12}$& $(\phi xy,r)$ & $ \eta_6$       &$\gamma_{12}$ &$(r\io)^{y\io}r^{xy\io}$ \\
$\alpha _{13}$& $(\phi yx,r)$ & $\eta_6$        &$\gamma_{13}$ & $(r\io)^{y\io}r^{x\io y\io}$\\
$\alpha _{14}$& $(\phi xy, s)$& $\eta_3 \eta_2$ &$\gamma_{14}$ &$(s\io)^{yxy\io}s^{xy\io}$ \\
$\alpha _{15}$& $(\phi xy,t)$ & $ \eta_1 \eta_2$&$\gamma_{15}$ & $ (t\io)^{y^{-2}}t^{xy\io}$\\
$\alpha _{16}$& $(\phi x^2,t)$& $\eta _5 $      &$\gamma_{16}$ & $(t\io)^{y\io}t^x$\\
$\alpha _{17}$& $(\phi yx,s)$ & $\eta_4 \eta_5$ &$\gamma_{17}$ & $ t^{x}(s\io)^x(t\io)^{y\io}s^{x\io y\io}$\\
$\alpha _{18}$& $(\phi yx,t)$ & $\eta_6 \eta_3\eta_7 \eta_5 $  & $\gamma_{18}$ &
           $t^x(s\io)^xr\io(s\io)^{x\io}ts\io(r\io)^{y\io}t^{x\io y\io}$\\
\hline
 \end{tabular}
\end{center}
We now let $I$ consist of elements $\ba_i$ in one-to-one
correspondence with the $\alpha_i$, and let $C_3(I)$ be the free
$G$-module on $I$. Define $\delta_3 \ba_i$ to be the value
$\gamma_i \in C_2$ given in the fourth column of the table. Let
$\wtI= G \times I$, and let $C_3(\wtI)$ be the free $\wtG$-module
on $\wtI$. We define $$ h_2(\alpha_i)=(1,\ba_i), i=1,\ldots,18. $$
Then we have $$h_1(\wt{\delta}_2 \alpha_i)\io \alpha_i^{h_0
\beta(\alpha_i)} = (1,\gamma_i) = \wt{\delta}_3 h_2 \alpha_i,
i=1,\ldots,18.$$ So we have extended our covering complex and its
contracting homotopy by one stage.

 However, we can in fact omit all of the
$\alpha_i$ except the first four, because of the trivial
identities $\gamma_5, \ldots,\gamma_{10}$, and the further
relations given in Table 2 below. We give the verification for the
last two further relations, the others being trivial or easy.

We note that
\begin{align*}
\gamma_{17}&= t^{x}(s\io)^x(t\io)^{y\io}s^{x\io y\io} \\
           &= t^{x}(t\io)^{y\io}(s\io)^{xyy\io x\io  y \io x\io y\io} s^{x\io y\io} \tag*{by the
crossed module rules} \\
            &= (t\io)^{y\io} t^{ y\io x\io y\io}((s\io)^{y\io}s)^{x \io
            y     \io}\\
           &= (t\io)^{y \io} t^{x} ((s\io)^{y\io}s)^{x \io y \io}
            \tag*{since $t^{(\delta_2t\io)x}=t^{x}$} \\
           &= \gamma_3   (\gamma_2\io)^{x\io y\io}
\end{align*}

In order to verify the further identity for $\gamma_{18}$, we
consider the abelianisation $C(R)^{\rm ab}$, which
is isomorphic to the free $S_3$-module on $R$. The difference
$\gamma_{18}-\gamma_4$ in $C(R)^{\rm ab}$  is
$$t.(\phi(x\io y\io)-\phi(x\io))= t.(\phi(xy)-1)\phi(x^2). $$
 Since the module of identities
is mapped injectively into $C(R)^{\rm ab}$, \cite{B-Hu}, and in $C(R)^{\rm ab}$ we have
$\gamma_3= t.(\phi(x)-\phi(y\io))$, the result follows.
So we have a set of four generators for
the module of identities for this presentation of $S_3$, of which
the first three given belong to the root module (see \cite{B-Hu} for
an account of this).

Let $J$ consist of the elements $\ba_i, i=1,\ldots,4$, and let $C_3(J)$
be the free $G$-module on $J$, with the restriction to it of the boundary $\delta_3$.
Let $\mathsf{r} : C_3(I) \to C_3(J)$ be be the $G$-module morphism defined by $\mathsf{r}(\ba_i)= \ba_i,
i = 1, \ldots, 4,\;  \mathsf{r}(\ba_i) =0, i =5, \ldots, 10$ and otherwise as
 in Table 2 below, so that $\delta_3 \mathsf{r}= \delta_3$. Note that $C_3(I)$ is an $S_3$-module, and so we write it
additively as a group, and use $.$ for the group action. To simplify the notation
we write these acting elements as words in  the generators $x,y$.

\begin{center}Table 2

\begin{tabular}{||l|l||l||}
\hline $\gamma_i$ & identity & $ \mathsf{r} \ba_i$ \\ \hline
    $\gamma_{11}$&= $(\gamma_1\io)^{x\io}$ & $ - \ba_1.
x^2$\rule{0em}{2.75ex}\\ $\gamma_{12} $&= $\gamma_1^{y\io}$ &
$\ba_1.y$\\ $\gamma_{13} $&= $(\gamma_1\io)^{x\io y\io}$
&$-\ba_1.yx$\\ $\gamma_{14} $&= $\gamma_2^{yxy\io}$ &$
\ba_2.x^2 $\\ $\gamma_{15} $&= $\gamma_3^{y\io}$ &$ \ba_3.y
$\\ $\gamma_{16} $&= $\gamma_3$ &$ \ba_3$\\ $\gamma_{17} $&=
$\gamma_3  (\gamma_2\io)^{x\io y\io}$ &$\ba_3-
\ba_2.yx$\\ $\gamma_{18} $&= $\gamma_4 \gamma_3^{yx^2}$ &$
\ba_4+ \ba_3.xy $ \rule{0em}{2.75ex}\\ \hline
\end{tabular}
\end{center}

 Define
$\mathsf{r}': C_3(\wtI) \to C_3(\wtJ), (g,d) \mapsto
(g,\mathsf{r}d)$, and define $h_2'= \mathsf{r}'h_2: C(\wtR) \to
C_3(\wtJ)$.
Then we have for $i=1,2,\ldots,18$
$$  \wt{\delta}_3 h_{2}' \A_i = h_1(\wt{\delta}_2 \A_i)\io\A_i^{h_0 \beta \A_i}
$$ and so we have a contracting homotopy up to this level.

Note that we now have $24$ generators $\tilde{d}=(g,\ba_i), g \in
S_3, i=1,\ldots,4$ of $C_3(\wtJ)$ and we can proceed to the next
stage, to obtain identities between identities corresponding to
each of these generators, namely $$p_3(-h_{2}\wt{\delta}_3
\tilde{d} +\tilde{d}^{h_0\beta \tilde{d}})$$ for
$\tilde{d}=(g,\ba_i), g \in S_3, i=1,\ldots,4$. This requires
another table.  In order to show how the calculations  go, we next
carry out  one intermediate calculation, and one full calculation.
 Further details of the calculations required for the table are omitted,
  but are available on request.

 Recall that $\kk$ is a groupoid morphism, and kills the
action  of $F(\wtX)$. So, for example,
\begin{alignat*}{2}
\kk (\phi y,t^x)&=\kk ((\phi yx\io,t)^{(\phi xy,x)})&& \\
             &= \kk (\phi xy,t)                  &&= \kk (\A_{15})\\
        &=(1,\ba_3^{y\io}) &&=(1,\ba_3^{\phi y\io})\\
        &= (1,\ba_3^y).&&
 \end{alignat*}
So we have
\Env{align*}{& -\kk \wt{\delta}_3 (\phi yx,\ba_4)+(\phi yx,\ba_4)^{(yx,x\io y\io)}\\
 &= -\kk (\phi yx,(s\io)^{x\io}ts\io (r\io)^{y\io}t^x(s\io)^xr\io t^{x\io})
+(1, \ba_4^{x\io y\io}) \\
&= -(1, -\ba_2.x^2 +\ba_4 +\ba_3.xy -\ba_3 +\ba_2.yx -\ba_1
  -\ba_2 +\ba_1.yx + \ba_3.y -  \ba_4.yx) \\
&= (1, \ba_1.(1-yx) +\ba_2.(1 +x^2- yx) +\ba_3.(1-y-xy) +\ba_4.(-1+yx))}

Some of the identities in Table 3 might seem as surprising to
others as they were to the authors. There is a process for
checking that these are identities among identities as follows.

We are required to check that $\delta_3$ of some combination $u$
of the $\ba_i$ is zero. Certainly each $\delta_3 \ba_i$ is an
identity among relations, and hence so is the corresponding linear
combination $u$. Therefore $u$ is 0 if and only if it maps to 0 in
the abelianised group $C(R)^{\mathrm{ab}}$, which is freely generated as a
${\mathbb Z}S_3$ module by the elements $r,s,t$. Thus we determine
the coefficients of these elements for the image of $u$ in
$C(R)^{\mathrm{ab}}$, and it is straightforward to check that these are
zero. This is analogous to a previous calculation.

\begin{center}
Table 3
\nopagebreak

\begin{tabular}{||l|l||l|l||}
\hline
\multicolumn{2}{||l||}{generator} &
\multicolumn{2}{c||}{identity $=p_3\left( -h_{2}\wt{\delta}_3 \theta _i
+\theta _i^{h_0\beta \theta _i}\right)$ \rule{0em}{3ex}}\\
\hline
$\xi _1$ & $(1,\ba_1)$ & $\mu_1$ & 0 \\
$\xi    _2$ & $(1,\ba_2)$ & $\mu_2$ & 0 \\
$\xi    _3$ & $(1,\ba_3)$ & $\mu_3$ & 0 \\
$\xi    _4$ & $(1,\ba_4)$ & $\mu_4$ & 0 \\
$\xi    _5$ & $(\phi x,\ba_1)$ & $\mu_5$ & 0 \\
$\xi    _6$ & $(\phi x,\ba_2)$ & $\mu_6$ &0  \\
$\xi    _7$ & $(\phi x,\ba_3)$ & $\mu_7$ &$\ba_3.(y+x^2)$  \\
$\xi    _8$ & $(\phi x,\ba_4)$ & $\mu_8$ &$\ba_1.(y-x^2)+\ba_4.(x^2-1)$  \\
$\xi    _9$ & $(\phi x^2,\ba_1)$ & $\mu_9$ &$\ba_1.(1+x+x^2)$  \\
$\xi    _{10}$ & $(\phi x^2,\ba_2)$ & $\mu_{10}$ &$\ba_2.(1+y)x$ \\
$\xi    _{11}$ & $(\phi x^2,\ba_3)$ & $\mu_{11}$ &$\ba_3.(1+yx)x$  \\
$\xi    _{12}$ & $(\phi x^2,\ba_4)$ & $\mu_{12}$ & $\ba_1.(1-yx) +\ba_4.(x-1)$ \\
$\xi    _{13}$ & $(\phi y,\ba_1)$ & $\mu_{13}$ &0  \\
$\xi    _{14}$ & $(\phi y,\ba_2)$ & $\mu_{14}$ &$\ba_2.(1+y)$  \\
$\xi_{15}$ & $(\phi y,\ba_3)$ & $\mu_{15}$ &0  \\
$\xi_{16}$ & $(\phi y,\ba_4)$ & $\mu_{16}$ &$ \ba_2.(1+x^2-yx)+\ba_3.(1-y -xy) -\ba_4.(1-y)$  \\
$\xi_{17}$ & $(\phi yx,\ba_1)$ & $\mu_{17}$ &0  \\
$\xi_{18}$ & $(\phi yx,\ba_2)$ & $\mu_{18}$ &0 \\
$\xi_{19}$ & $(\phi yx,\ba_3)$ & $\mu_{19}$ & $\ba_3.(1+yx)$ \\
$\xi_{20}$ & $(\phi yx,\ba_4)$ & $\mu_{20}$ &$\ba_1.(1-yx) +\ba_2.(1 +x^2- yx) +\ba_3.(1-y-xy) +\ba_4.(-1+yx)$  \\
$\xi_{21}$ & $(\phi yx^2,\ba_1)$ & $\mu_{21}$ &$\ba_1.(1+x+x^2)y$ \\
$\xi_{22}$ & $(\phi yx^2,\ba_2)$ & $\mu_{22}$ & $\ba_2.(1+y)x^2$ \\
$\xi_{23}$ & $(\phi yx^2,\ba_3)$ & $\mu_{23}$ & 0 \\
$\xi_{24}$ & $(\phi yx^2,\ba_4)$ & $\mu_{24}$ &$ \ba_1.(y- x^2) +\ba_2.(1+x^2-yx)
                                     +\ba_3.(1+x^2-xy) +\ba_4.(-1+xy)$ \\

\hline

\end{tabular}

\end{center}

We next reduce this to a smaller, and clearly minimal,   set of identities
among identities, as in the following table.

\begin{center}Table 4
\nopagebreak

\begin{tabular}{||l|l||}
\hline
generator & definition/further identity \\
\hline
$\mu_9   $  &  $\ba_1.(1+x+x^2)$ \\
$\mu_{14}$  &   $\ba_2.(1+y)$\\
$\mu_{19}$  &   $\ba_3.(1+yx)$ \\
$\mu_{12}$  &   $\ba_1.(1-yx) +\ba_4.(x-1)$\\
$\mu_{16}$  &  $\ba_2.(1+x^2-yx)+ \ba_3.(1-y -xy) -\ba_4.(1-y)$  \\
$\mu_{21}$  & $=\mu_{9}.y$  \\
$\mu_{10}$  & $=\mu_{14}.x$  \\
$\mu_{22}$  & $=\mu_{14}.x^2$   \\
$\mu_{11}$  & $=\mu_{19}.x$  \\
$\mu_7   $  & $=\mu_{19}.x^2$  \\
$\mu_8   $  &  $=\mu_9.(-1+y)+\mu_{12}.(x+1)$ \\
$\mu_{20}$  &  $=\mu_9.(1-y) +\mu_{12}.(x+1)y +\mu_{16}$ \\
$\mu_{24}$  &  $= \mu_{12}.y +\mu_{16}+\mu _{19} . x^2$ \\

\hline
\end{tabular}
\end{center}
Let $K$ be set with elements $\bar{\mu}_{9},
\bar{\mu}_{14},\bar{\mu}_{19}$, $ \bar{\mu}_{12}, \bar{\mu}_{16}$,
let $C_4(K)$ be the free $G$-module on $K$, and let $\delta_4:
C_4(K) \to C_3(J)$ be given by the first four line of the second
column of Table 4. Then the sequence $C_4(K) \to C_3(J)\to C(R)$
is exact and we have extended our crossed resolution by one
further step. Hence we have a presentation of the $G$-module
$\pi(\cP)$.

Such a crossed resolution has been extended by two further steps,
but with different choices, in \cite{He}.

As explained in the Introduction, this example is chosen as  one
which illustrates the method, which has non trivial calculations
but also is perhaps the largest example of this type which
one would care to do by hand. The major problems are the
calculation of $h_1$, i.e. representing a set of group generators of
$N=\Ker \phi$ as consequences of the relators, and more seriously,
calculating minimal generating subsets of sets of generators of
submodules of free $\Z G$-modules,  as
well as finding the relations giving all the generators in terms of the
smaller set. The first problem is dealt with in \cite{HeW} and the
second in \cite{HeR}.

\section{Presentations of groupoids}
\label{presgpds}
The category of groupoids will  be written ${\Gpds}$.
Our convention for groupoids is that the composite of arrows $a:
x\to y, b: y \to z$ is written $ab: x \to z$.

The theory of groupoids may be thought of as an algebraic
analogue of the theory of groups, but based on directed graphs
rather than on sets. For some discussion of the philosophy of this,
see \cite{B-Irish}.
\subsection{Free groupoids}
The term {\em graph} will always mean what is usually called a
directed graph.
A {\em  graph} $X$ consists of two sets
$Arr(X),Ob(X)$, of arrows and objects respectively of
$X$, and two functions $s,t : Arr(X) \to Ob(X)$,
called the {\em source} and {\em target} maps. A {\em morphism}
$f : X \to Y$ of  graphs consists of two functions $Arr(X)
\to Arr(Y), Ob(X) \to Ob(Y)$, which commute with the source and
target maps. This defines the category $\DirG$.

A basic construction in any algebraic theory is that of free
objects. For groups, the free group functor $F: \Sets \to \Groups$
is left adjoint to the forgetful functor $\Groups \to \Sets$. In
the case of groupoids, we may define the {\em free groupoid}
functor to be the left adjoint $F: \DirG\to \Gpds$ to the
forgetful functor $U:\Gpds \to \DirG$ giving the underlying graph
$UG$ of a groupoid $G$, namely forgetting the composition, the identity function
$Ob(G) \to G$, and the inverse map $G \to G$.
 So if $X$ is a graph, then the free
groupoid $F(X)$ on $X$ consists of a graph morphism $i: X \to
UF(X)$ which is universal for morphisms from $X$ to the underlying
graph of a groupoid.

The set of objects of $F(X)$ may be identified with $Ob(X)$. There
are several ways of explicitly constructing the set of arrows of
$F(X)$.  The usual way is as equivalence classes of {\em composable} words
$$w=(x_1,\eps_1)\ldots (x_n,\eps_n), n \ge 0,x_i \in Arr(X), \eps=\pm $$
together with empty words $(\,)_a, a \in Ob(X)$, where the word $w$
is composable means that $t(x_i,\eps_i) = s(x_{i+1}, \eps_{i+1}),
i=1\ldots n-1$, where
$$s(x,\eps)=\begin{cases} sx &\text{if } \eps =+,\\
                          tx &\text{if } \eps =-,\end{cases}
\qquad t(x,\eps)=\begin{cases} tx &\text{if } \eps =+,\\
                          sx &\text{if } \eps =-.\end{cases}   $$
The equivalence relation on words, and the composition, to obtain
the free groupoid is defined in a manner analogous to the usual
definition of free group, and the graph morphism $i: X \to F(X)$
sends an arrow $x$ to $[x]$, the equivalence class of the word
$(x,+)$.

A groupoid $G$ is called {\em connected} if $G(a,b)$ is non empty
for all $a,b \in Ob(G)$. The maximal connected subgroupoids of $G$
are called the {\em (connected) components} of $G$.

If $a$ is an object of the groupoid $G$, then the set $G(a,a)$
inherits a group structure from the composition on $G$, and this is
called the {\em object group} of $G$ at $a$ and is written also
$G(a)$. The groupoid $G$ is called {\em simply connected} if all its
object groups are trivial. If it is connected and simply connected,
it is called 1-{\em connected}, or a {\em tree groupoid}.

A standard example of a tree groupoid is the {\em indiscrete}, or
{\em square}, groupoid $I(S)$ on a set $S$. This has object set
$S$ and arrow set $S \times S$, with $s,t :S\times S \to S$ being
the first and second projections. The composition on $I(S)$ is given
by $$(a,b)(b,c) = (a,c), a,b,c \in S. $$

A  graph $X$ is called {\em connected} if the free groupoid $F(X)$
on $X$ is connected, and is called a {\em forest} if every object
group $F(X)(a)$ of $F(X), a \in Ob(X),$ is trivial. A connected
forest is called a {\em tree}. If $X$ is a tree, then  $F(X)$ is a
tree groupoid.

\subsection{Retractions}
Let $G$ be a connected groupoid. Let $a_0$ be an object of $G$.
For each object $a$ of $G$ choose an arrow $\tau a : a \to a_0$,
with $\tau a_0 =1_{a_0}$. Then an isomorphism $$\phi : G \to
G(a_0) \times I(Ob(G))$$ is given by $g \mapsto ((\tau a)\io g (\tau
b), (a,b)), g \in G(a,b), a,b \in Ob(G). $ The composition of
$\phi$ with the projection yields a morphism $\rho : G \to G(a_0)$
which we call a {\em deformation retraction}, since it is the
identity on $G(a_0)$ and is in fact homotopic to the identity
morphism of $G$, though we do not elaborate on this fact here.

It is also standard \cite[8.1.5]{B-88} that a connected groupoid
$G$ is isomorphic to the free product groupoid $G(a_0) \ast T$
where $a_0 \in Ob(G)$ and $T$ is any wide, tree subgroupoid of
$G$. The importance of this is as follows.

Suppose that $X$ is a graph which generates the connected groupoid
$G$. Then $X$ is connected. Choose a maximal tree $T$ in $X$. Then
$T$ determines for each $a_0$ in $Ob(G)$ a retraction $\rho_T: G
\to G(a_0)$ and the isomorphisms $$G \cong G(a_0) \ast
I(Ob(G))\cong G(a_0) \ast F(T)$$ show that a morphism $G \to K$
from $G$ to a groupoid $K$ is completely determined by a morphism
of groupoids $G(a_0) \to K $ and a graph morphism $T \to K$ which
agree on the object $a_0$.

We shall use later the following proposition,
which is a special case of   \cite[6.7.3]{B-88}:
\begin{prop} \label{pushout} Let $G,H$ be
groupoids with the same set of objects, and let
$\phi: G \to H$ be a morphism of groupoids which is the
identity on objects. Suppose that $G$ is connected and $a_0 \in
Ob(G)$. Choose a retraction $\rho : G \to G(a_0)$. Then there is a
retraction $\sigma : H \to H(a_0)$ such that the following
diagram, where $\phi'$ is the restriction of $\phi$:
\begin{equation}
{\sqdiagram{G}{\rho}{G(a_0)}{\phi}{\phi'}{H}{\sigma}{H(a_0)}}
\end{equation}
is commutative and is a pushout of groupoids.
\end{prop}

\subsection{Normal subgroupoids and quotient groupoids}
Let $G$ be a groupoid.  A subgroupoid $N$ of $G$ is called {\em
normal} if $N$ is wide in $G$ (i.e. $Ob(N) = Ob(G)$) and for
any objects $a,b$ of $G$ and
$g$ in $G(b,a), \; g^{-1}N(b)g =  N(a)$.

Let $\phi : G \to H$ be a morphism of groupoids.  Then $\Ker \phi$
is the wide subgroupoid of $G$ whose elements are all $g$ in $G$
such that $\phi g$ is an identity of $H$ is a normal subgroupoid
of $G$. If $Ob (f)$ is injective then $\Ker \phi$ is totally
disconnected, i.e. $(\Ker \phi) (a,b)= \emptyset $ if $a \ne b$.

A morphism $\phi : G \to H$ is said to {\em annihilate} a subgraph
$X$ of $G$ if $\phi(X) $ is a discrete subgroupoid of $H$.  Thus
$\Ker \phi$ is the largest subgroupoid of $G$ annihilated by
$\phi$. The next proposition gives the existence of quotient
groupoids.
\begin{prop}
Let $N$ be a totally disconnected, normal subgroupoid of $G$. Then
there is a groupoid $G/N$ and a morphism $p : G \to G/N$ such that
$p$ annihilates $N$ and is universal for morphisms from $G$ which
annihilate $N$.
\end{prop}

\begin{proof}
We define $Ob (G / N)= Ob (G)$. If $a,b \in Ob (G)$ we define
$(G/N)(a,b)$ to consist of all cosets $N(a)g, g \in G(a,b). $ The
multiplication of $G$ is inherited by $G/N$, which becomes a
groupoid.

The morphism $p: G \to G/N$ is the identity on objects, and on
elements is defined by $g \mapsto N(sg)g$. Clearly $p$ is a
morphism and $\Ker p = N$.

The remainder of the proof is clear.
\end{proof}

We call $G/N$ a {\em quotient groupoid} of $G$.

\subsection{Presentations of groupoids}
We now consider relations in a groupoid. Suppose given for each
object $a$ of the groupoid $G$ a set $R(a)$ of elements of
$G(a)$--- thus $R$ can be regarded as a wide, totally disconnected
subgraph of $G$. The {\em normal closure} $N(R)$ of $R$ is the
smallest wide normal subgroupoid of $G$ which contains $R$. This
obviously exists since the intersection of any family of normal
subgroupoids of $G$ is again a normal subgroupoid of $G$. Further,
$N(R)$ is totally disconnected since the family of object groups
of any normal subgroupoid $N$ of $G$ is again a normal subgroupoid
of $G$.

Alternatively, $N = N(R)$ can be constructed explicitly. Let $a$
be an object of $G$. By a {\em consequence} of $R$ at $a$ is meant
either the identity of $G$ at $a$, or any product
\begin{equation} \label{conseq}
    \tau= g_{1}^{-1}r_{1}^{\eps_1}g_1 \ldots
g_{n}^{-1} r_{n}^{\eps_n} g_{n},
\end{equation}
in which $ n \ge 1,\;g_{i} \in G(a_{i},a)$ for some object $a_i$
of $G, \; \eps_i= \pm 1$ and $r_{i}$  is an element of $R(a_i)$.
Clearly, the set $N(a)$ of consequences of $R$ at $a$ is a
subgroup of $G(a)$ and the family $N= (N(a):a \in Ob(G))$ of these
groups is a totally disconnected normal subgroupoid of $G$
containing $R$. Clearly $N=N(R)$.

The projection $p:G \to G/N(R)$ has the following universal
property: {\em if $f : G \to H$ is any morphism which annihilates
$R$ then there is a unique morphism $f
: G / N(R)
\to H$ such that $f p = f$.}  We call $G / N(R)$ {\em the groupoid
$G$ with the relations $ r= 1, r\in R$.}

In applications, we are often given $G$, $R$ as above and wish to
describe the object groups of $G/N(R)$. These are determined by
the following result. \begin {prop}{\em \cite[8.3.3]{B-88}} Let G
be connected, let $a_0 \in Ob(G)$ and let $\rho :G \to G(a_0)$ be
a deformation retraction. Let $H = G/N(R)$. Then $H(a_0)$ is
isomorphic to the group $G(a_0)$ with the relations
\[
\rho (  r) = 1 ,    r \in R.
\]
\end {prop}
\begin{proof}
The proof follows  from Proposition \ref{pushout}, with
$H=G/N$ and $\phi=p:G\to G/N$ the quotient morphism. Details are given
in \cite{B-88}.
\end{proof}

\section{Crossed modules and free crossed modules\\
 over groupoids}
The theory of crossed modules and free crossed modules is due to
Whitehead \cite{W1}. Expositions are given in for example
\cite{B-Hu,HMS}. In order to obtain an algebraic model of
universal covers,  we need the corresponding definitions for the
groupoid case,  due to Brown and Higgins in \cite{BH2}.

Let $\Phi$ be a groupoid. A {\em crossed $\Phi$-module}
consists of: \begin{enumerate}[(i)]
\item   a
totally disconnected groupoid $M$ with the same object set as $\Phi$;
\item a morphism $\mu : M \to \Phi$ of groupoids which is the
identity on objects; and
\item an action of the groupoid $\Phi$ on the right of the
groupoid $M$ {\em via} $\mu$.  \end{enumerate}
This last condition means that if $x \in \Phi(a,b), m \in M(p)$, then $m^x \in
M(b)$ and the usual laws of an action apply, namely  $m^1=m,
(m^x)^y=m^{xy}, (mn)^x=m^xn^x$ whenever the terms are defined.

The axioms for a crossed module are:\begin{enumerate}[CM1)]
\item  $\mu (m^x) = x \io (\mu m)x$,
\item  $ n\io m
n = m^{\mu n}$,
 \end{enumerate}
for all $m,n \in M, x \in \Phi$ and whenever the terms are defined.

Such a crossed $\Phi$-module is written $(M,\mu, \Phi)$ or
$\mu : M \to \Phi$, or simply as $M$.

A {\em morphism} from a crossed module $\mu : M\to \Phi$ to a crossed
module  $\nu : N \to \Psi$ consists of a pair of morphisms of
groupoids $f: \Phi \to \Psi, \; g : M \to N$ such that $\nu g = f
\mu$ and $g(m^x) = (gm)^{fx}$ whenever $m^x$ is defined. This yields
the category ${\XMod}$ of crossed modules and their morphisms.

There is also a category ${\PXMod}$ of precrossed modules, in which the
axiom CM2) is dropped. The inclusion of categories ${\XMod} \to {\PXMod}$
has a left adjoint constructed as follows.

Let $\mu : M \to \Phi$ be a precrossed module. By a {\em Peiffer element},
or {\em twisted commutator}, is meant an element
$$\langle m,n\rangle = m\io n\io m n^{\mu m} $$
where $m,n\in M(p)$ for some object $p$. As in the group case
(see \cite[Proposition 2, p.158]{B-Hu}) one proves that the Peiffer elements
generate a normal $\Phi$-invariant subgroupoid $\langle M,M\rangle$ of $M$,
and the quotient groupoid, $M ^{ass}=M /\langle M,M \rangle$, with
the induced morphism $\mu ' : M ^{ass} \to \Phi$, inherits the
structure of crossed module. This {\em associated crossed module}
gives the reflection from the category ${\PXMod}$ of
precrossed modules to the category ${\XMod}$ of crossed modules as
required.

Let $\Phi$ be a groupoid, let $R$ be a totally disconnected graph
with the same object
set as $\Phi$, and let $w: R \to \Phi$ be a graph morphism which is
the identity on objects. We define the {\em free crossed module on}
$w$ to be a crossed module $\partial : C(w) \to \Phi$
together with a graph morphism $\bar{w}: R \to C(w)$ such that:
\begin{enumerate}[(i)]
\item  $\partial \bar{w} = w$;
\item if $\mu : M \to \Phi$ is a crossed module and $g:R \to M$ is a
graph morphism over the identity on objects such that $\mu g =
w$, then there is a unique morphism $g' : C(R) \to M$ of crossed
$\Phi$-modules such that $g' \bar{w} = g$.
\end{enumerate}

Free crossed modules over groups were defined and constructed by
Whitehead \cite{W1}, and an exposition is given  in \cite{B-Hu}.
The analogous construction  for   groupoids  is as follows.

Let $w: R \to \Phi$ be given as above. One first forms
 the free groupoid $H(w)$ on the totally
disconnected graph $Y$ with object set $Ob(\Phi)$ where $Y(p)$
consists of pairs $(r,u)$ such that $r \in R(q), u \in \Phi(q,p)$.
Let $\partial':H(w) \to \Phi$ be given by $(r,u) \mapsto u\io (wr) u$,
and let $\Phi$ operate on $H(w)$ by $(r,u) ^v = (r,uv)$. This yields  the
{\em free precrossed module} on $w$, and the free crossed module is
the associated crossed module  $\partial: C(w)= H(w)^{ass}\to \Phi$.

Notice that the image $\partial (C(w))$ is the normal closure
$N(wR)$ of $wR$ in $\Phi$.

It is useful to see this construction as a special case of the {\em
induced crossed modules} of \cite{BH1} (but for the groupoid case),
which can be regarded
as arising from a `change of base' \cite{B-94}.
That is, $C(w)$ is isomorphic to the crossed module
$\omega_{\ast}F(R)$ induced from the identity crossed module $1: F(R)
\to F(R)$ by the morphism $\omega : F(R) \to \Phi$
determined by $w :R \to \Phi$. Further, we have a pushout of crossed modules
$$\sqdiagram{(1,0,F(R))}{(1,\omega)}{(1,0,\Phi)}{}{}{(F(R),1,F(R))}{}{(C(w), \partial,\Phi)}$$
This allows a link with the 2-dimensional Van Kampen Theorem
of \cite{BH1} (or rather, with the groupoid version formulated in all dimensions in
\cite{BH3}), to obtain a proof of a groupoid version of a well known theorem of
Whitehead \cite{W1}, as follows:
\begin{thm} Let $U_0$ be  a subset of the space $U$ and
suppose the space $V$ is obtained from $U$ by attaching 2-cells by maps
of pairs $f_r: (S^1,1) \to (U,U_0), r \in R$. Then the family
of second relative homotopy
groups $\pi_2(V,U,p), p\in U_0$ form  the free crossed  module
over the fundamental groupoid $\pi_1(U,U_0)$ on the graph morphism
$w: R \to \pi_1(U,U_0)$ given by $wr= (f_r) _{\ast}(\iota)$, where  $\iota$
here denotes a generator of the fundamental group $\pi_1(S^1,1)$.
\end{thm}

\section{Crossed complexes} \label{crscom}
The basic geometric example of a crossed complex  is the
{\em fundamental crossed
complex }
$\pi X\sast $ of a filtered space
$$  X \sast : X_0 \subseteq X_1 \subseteq  \cdots  \subseteq X_n
\subseteq   \cdots  \subseteq X.$$
Here $\pi_1X \sast $ is the fundamental
groupoid $\pi_1(X_1,X_0)$ and for $n\geq2$, $\pi_n X \sast $
  is the family of relative
homotopy groups $\pi_n(X_n,X_{n-1},p)$ for all $p \in X_0$. These
come equipped
with the standard operations of $\pi_1 X \sast $ on $\pi_n X \sast $
and boundary
maps $\delta : \pi_n X \sast  \to \pi_{n-1}X \sast $. The axioms for
crossed complexes
are those universally satisfied for this example.

The definition of a crossed complex generalises to the
case of a set of base points definitions given by Blakers
\cite{Blakers} (under the term `group system') and Whitehead \cite{W1},
under the term `homotopy system' (except that he restricted also to
the free case). We recall this general definition from
\cite{BH2}.

A {\em crossed complex } $C$ (over a groupoid) is a sequence of
morphisms of groupoids over $C_0$
$$\diagram
\cdots \rto & C_n \dto<-.05ex>^{\beta}\rto^{\delta_n} & C_{n-1} \rto \dto<-1.2ex>^{\beta} & \cdots \rto &
 C_2\rto^{\delta_2} \dto & C_1
\dto<0.0ex>^(0.45){\delta^1} \dto<-1ex>_(0.45){\delta^0} \\ & C_0&C_0\rule{0.5em}{0ex}  &
& \rule{0.5em}{0ex} C_0 &  \rule{0em}{0ex} C_0 . \enddiagram$$
Here $\{ C_n\}_{n \ge 2}$ is a family of groups with base point map $\beta$,
 and $ \delta^0, \delta^1$ are the source and targets
for the groupoid $C_1$. We further require given an operation of the groupoid
$C_1$ on each family of groups $C_n$ for
$n \ge 2$ such that:
\renewcommand{\labelenumi}{(\roman{enumi})}
   \begin{enumerate}
      \item  each $\delta_n$ is a morphism over the identity on $C_0$;

      \item  $C_2 \rightarrow C_1$ is a crossed module over $C_1$;

      \item  \label{modul} $C_n$ is a $C_1$-module for $n\geq3$;

      \item  $\delta : C_n \rightarrow C_{n-1}$ is an operator morphism
 for $n\geq3$;

      \item  $\delta\delta :C_n \rightarrow C_{n-2}$ is trivial for
$n\geq3$;

      \item  $\delta C_2$ acts trivially on $C_n$ for $n\geq3$.
   \end{enumerate}
Because of axiom (iii) we shall write the composition in
$C_n$ additively for $n \ge 3$, but we will use multiplicative notation
in dimensions 1 and 2.

Let $C$ be a crossed complex. Its {\em  fundamental groupoid} $\pi_1{C}$ is
the quotient of the
groupoid $C_1$ by the normal, totally disconnected subgroupoid
$\delta{C_2}$. The rules
for a crossed complex  give $C_n$, for $n\geq 3$, the induced structure
 of $\pi_1{C}$-module.

A {\it morphism} $f: C \rightarrow D$ of crossed complexes is a family
 of groupoid
morphisms $f_n : C_n \rightarrow D_n ~(n\geq 0)$ which preserves all the
structure. This defines the category ${\Crs}$ of crossed complexes.
The fundamental groupoid now gives a functor $ \pi_1 : {\Crs}\to {\Gpds}$. This functor
 is left adjoint to the functor $i : {\Gpds} \to {\Crs} $ where for a groupoid $G$ the crossed complex  $iG$
agrees with $G$ in dimensions 0 and 1, and is otherwise trivial.

An $m$-{\em truncated} crossed complex $C$ consists of all the structure defined above
but only for $n \le m$. In particular, an $m$-truncated crossed complex is for
$m=0,1,2$ simply a set, a groupoid, and a crossed module respectively.

\section{Covering morphisms of groupoids and crossed complexes} \label{cov}
For the convenience of readers, and to fix the notation, we recall
here the basic facts on covering morphisms of groupoids.

Let $G$ be a groupoid. For each object $a$ of $G$ the {\em star} of $a$ in $G$,
denoted by ${\rm St}_{G}~a$, is the union of the sets $G(a,b)$ for all objects $b$
of $G$, i.e. ${\rm St}_{G}a = \{g \in G : sg = a\}$.
A morphism $p : \wt{G} \rightarrow G$ of groupoids is a
{\em covering morphism} if
for each object $\wt{a}$ of $\wt{G}$ the restriction of $p$
\[{\rm St}_{\wt{G}}~\wt{a} \to {\rm St}_{G}~p\wt{a} \]
is bijective. In this case $\wt{G}$ is called a
{\em covering groupoid of G}.

A basic result for covering groupoids is { \em unique path lifting}.
That is, let $p : \wt{G} \rightarrow G$ be a covering morphism of groupoids,
and let $(g_1, g_2, \ldots, g_n)$ be a sequence of composable elements of $G$.
Let $\tilde{a} \in Ob(\wt{G})$ be such that $p\tilde{a} $ is the starting point of
$g_1$. Then there is a unique composable sequence
 $(\tilde{g}_1, \tilde{g}_2, \ldots, \tilde{g}_n) $ of elements of $\wt{G}$ such that
$\tilde{g}_1$ starts at $\tilde{a}$ and $p\tilde{g}_i = g_i, i=1, \ldots, n$.

If $G$ is a groupoid, the category ${\Gpds\Cov}/G$ of coverings
of $G$ has as
objects the covering morphisms $p : H \rightarrow G$ and has as arrows
(morphisms) the commutative diagrams of morphisms of groupoids, where $p$ and $q$ are
covering morphisms,
$$\diagram
  H \drto_p  \rrto^f  &         & K \dlto^q  \\
                      & G &
\enddiagram$$
By a result of \cite{B-88}, $f$ also is a covering morphism. It is
convenient to write such a diagram as a triple $(f,p,q)$. The composition
in ${\Gpds\Cov}/G$ is then given as usual by
\[(g,q,r)(f,p,q) = (gf,p,r). \]

It is a standard result (see for example \cite{Hi,B-70}) that
 the category
 ${\Gpds\Cov}/G$ is
equivalent to the functor category ${\Sets}^G.$ This is useful for
constructing covering morphisms of the groupoid $G$. For example,
if $a$ is an object of the transitive groupoid $G$, then the groupoid
$G$ operates on the family of stars ${\rm St}_G~a $, and the
associated covering morphism $\wt{G} \to G $ defines the {\em universal
cover} $\wt{G}$ of the groupoid $G$. In particular, this gives the
universal covering groupoid of a group.

We now give the generalisation of  this notion to crossed complexes.

\begin{Def}{\em   \cite{Howie} A  morphism $p : \wt{C} \rightarrow C$
 of crossed
complexes is a {\it covering morphism} if
\renewcommand{\theenumi}{(\roman{enumi})}
   \begin{enumerate}
      \item
 the morphism $p_1 : (\wt{C}_1,\wt{C}_0) \to (C_1,C_0)$ is a
covering morphism of groupoids;
       \item
 for each $n\geq 2$ and $\wt{x} \in \wt{C}_0$, the morphism of groups
$p _n : \wt{C}_n(\wt{x}) \rightarrow C_n
(p\wt{x})$ is an isomorphism.
  \end{enumerate}}
\end{Def}
In such case we call $\wt{C}$ a {\it covering crossed complex } of $C$.

This definition may also be expressed in terms of the unique covering
homotopy property. For more details (but there with emphasis on
fibrations) see \cite{BH-91}.

\Env{prop}{\label{isocov} Let $p: \wt{C} \to C$ be a covering morphism of crossed
complexes and let $\tilde{a}\in Ob(\wt{C}). $ Let $a=p\tilde{a}, $
and let $K=p_0\io (a) \subseteq Ob( \wt{C}).$
Then $p$ induces isomorphisms $\pi_n(\wt{C},\tilde{a})\to \pi_n(C,a)$
for $n \ge 2$ and a sequence
$$1 \to \pi_1(\wt{C}, \tilde{a}) \to \pi_1(C,a) \to K \to \pi_0(\wt{C}) \to \pi_0(C) $$
which is exact in the sense of the exact sequence of a fibration of groupoids.
}
The comment about exactness has to do with operations on the pointed sets: see
\cite{B-70,B-88}. The proof of the proposition is easy and is omitted.

The following result gives a basic geometric example of
a covering morphism of crossed complexes.
\begin{thm} Let $X \sast$  and  $ Y\sast $  be filtered spaces
and let $f : X \to Y$ be a
covering map of spaces such that for each $n \geq 0$, $f_n : X_n \to Y_n$
is also
a covering map with $X_n = f^{-1}(Y_n)$. Then $\pi f : \pi X \sast  \to \pi Y \sast$ is a
covering morphism of crossed complexes.
\end{thm}

\begin{proof} By a result of \cite{B-88}, $\pi f:\pi_1 X_1 \to \pi_1 Y_1$
 is a covering
morphism of groupoids. Since $X_0 = f^{-1}(Y_0)$, the restriction
 of $\pi_1 f $ to $\pi_1(X_1,X_0)
\to \pi_1(Y_1,Y_0)$ is also a covering morphism of groupoids.
Now for each $n \geq 2$
and for each $x_0 \in X_0$, $f\sast : \pi_n(X_n,X_{n-1},x_0) \to
\pi_n(Y_n,Y_{n-1},p(x_0))$
is an isomorphism (see for example, \cite{Hu}).
\end{proof}

Here is an important method of constructing new covering  morphisms.

\begin{prop} Let $p : \wt{C} \to C$ be a covering
 morphism of crossed complexes. Then
the induced morphism $\pi_{1}(p): \pi_{1}\wt{C} \to \pi_{1}C$
is a covering morphism
of groupoids. \label{cov-c-g}
\end{prop}

\begin{proof}
Let $ \tilde{x} \in \wt{C}_0$. We will show that $p_{\tilde{x}}^{'}  :
{\rm St}_{\pi_1 \wt{C}}~ \tilde{x} \rightarrow
{\rm St}_{\pi_1C}~ p\tilde{x}$
is bijective. Let $[a] \in {\rm St}_{\pi_1C}~p\tilde{x}$,
where  $a \in {\rm St}_{C}
p\tilde{x}$.
Since $p$ is a covering morphism, there exists a unique $\tilde{a}$
of
${\rm St}_{\wt{C}}{\tilde{x}}$ such that $p \tilde{a} = a$. So
$p_{\tilde{x}}^{'}[\tilde{a}] = [a]$ and thus $p_{\tilde{x}}^{'}$ is
surjective.

Now suppose that
 $p_{\tilde{x}}^{'}[\tilde{a}] = p_{\tilde{x}}^{'}[\tilde{b}]$.
Then
$(p\tilde{b})\io  p\tilde{a} \in \delta C_2(p\tilde{x})$
which implies that $(p\tilde{b})\io( p\tilde{a}) = \delta  p\tilde{c}$
for a unique $\tilde{c} \in \wt{C}_2(\tilde{x})$.
Because $p$ is a covering morphism, we need only show that $(\tilde{b})\io
\tilde{a} = \delta\tilde{c}$. This follows by star injectivity.
Therefore $p_{\tilde{x}}^{'}$ is injective and so is bijective.
 Hence $\pi_{1}(p) $
is a covering morphism of groupoids.
\end{proof}

Let $C$ be a crossed complex. We write  ${\Crs\Cov}/C$ for
the full subcategory of the slice category ${\Crs}/C $
whose objects are the covering morphisms of $C$.

\begin{prop} \label{pullcov} Suppose given a pullback diagram
of crossed complexes
$$  \sqdiagram{\wt{C}}{\bar{f}}{\wt{E}}{\bar{q}}{q}{C}{f}{E}
$$
in which $q$ is a covering morphism. Then
$\bar{q}$ is a  covering morphism.
\end{prop}

We omit the proof. The groupoid case is done in
 \cite[9.7.6]{B-88}.  See also \cite{BHK} for uses of
pullbacks of covering morphisms of groupoids.

Our next result is the analogue for covering morphisms of
crossed complexes  of a
classical  result for covering maps of spaces \cite[9.6.1]{B-88}.

\begin{thm} \label{equiv} If $C$ is a crossed complex, then the functor
$\pi_1 : {\Crs}\rightarrow {\Gpds}$ induces an equivalence of categories
  $$\pi_1^{'} : {\Crs\Cov}/C \rightarrow {\Gpds\Cov}/(\pi_1C).$$
\end{thm}

\begin{proof} If $p : \wt{C} \to C$  is a covering morphism
 of crossed complexes, then
$\pi_{1}p : \pi_{1}{\wt{C}} \to \pi_1{C}$ is a covering
 morphism of groupoids, by
Proposition \ref{cov-c-g}.
Since $\pi_{1}$ is a functor, we also obtain the functor
$\pi_{1}^{'}$. To prove $\pi_{1}
^{'}$ is an equivalence of categories, we construct a
functor $\rho : {\Gpds\Cov}/
(\pi_{1}C) \to {\Crs\Cov}/C$ and prove that there
are equivalences of functors  $1
\simeq  \rho\pi_{1}^{'}$ and $1 \simeq \pi_{1}^{'}\rho.$

Let $C$ be a crossed complex, and let $q : D \to \pi_1C$ be a
covering morphism of groupoids. Let $\wt{C}$ be given by the
 pullback diagram in
the category of crossed complexes:
\begin{equation}{
 \sqdiagram{\wt{C}}{\bar{\phi}}{  iD}{\bar{q}}{q}{C }{\phi}{i\pi_1 C}}
 \label{pullcov2} \end{equation}
By proposition \ref{pullcov}, $\bar{q} : \wt{C} \to C$ is a
 covering morphism of  crossed complexes.

We define  the functor $\rho$ by $\rho (q)  = \bar{q}  $, and
extend $\rho $ in the obvious way to morphisms.

 The natural transformation
$\pi_1 ' \rho \simeq 1$ is defined on a covering morphism
$q : D \to \pi_1C$ to be the composite morphism
 $$\lambda: \pi_1(\wt{C})\llabto{1}{\pi_1(\bar{\phi})} \pi_1(iD)\cong D$$ where
$\bar{\phi}:\wt{C}\to iD $ is given in  diagram
(\ref{pullcov2}).
The proof  that $\lambda$ is an isomorphism is simple and is left
to the reader.

To prove that $1 \simeq \rho \pi_1'$,  we show that the
following diagram is a pullback:
$$ \sqdiagram{\wt{C}}{\wt{\phi}}{i\pi_1\wt{C}}{q}{q'=i\pi_1(q)}{C}{\phi}{i\pi_1C} $$
This is clear in dimension 0 and in dimensions
$\ge 2 $. For the case of dimension 1, let $c : x
\to y $ in $C$, and $[\tilde c] \in (\pi_1\wt{C})
(\tilde x,\tilde y)$ be such that $q[\tilde
c]=\phi(c)$. Then there exists a unique $\tilde c'
: \tilde x \to \tilde y$ such that
$\wt{\phi}(\tilde c')=[\tilde c]$ and $\bar{q}(\tilde c') = c.$
Now, $\bar{q}(\tilde{c}\delta \wt C_2(\tilde{x})) = \phi(c) = c
\delta C_2(x).$  This implies that $(\bar{q}\tilde{c}) \delta
C_2(x)=c \delta C_2(x)$. So $\bar{q}(\tilde c)=c(\delta c_2)$ for
some $c_2 \in C_2(x).$ Therefore there exists a unique $\tilde c_2
\in \wt C_2(\tilde x)$ covering $c_2$, and
$\bar{q}(\tilde{c}(\delta \tilde c_2)^{-1})=c.$ So the above
diagram is a pullback and thus we have proved that $1 \simeq \rho
\pi_1'$. This proves the equivalence of the two categories.
\end{proof}

\section{Covering morphisms and colimits } \label{adj}

In this section we give a result due to Howie \cite[Theorem
5.1]{Howie} which  we use to prove covering crossed complexes of
free crossed complexes are free.

\begin{thm}
Let $p: A \to B$ be a morphism of crossed complexes. Then $p$ is a
fibration if and only if the pullback functor $p^{\ast}: {\Crs}/B
\to {\Crs}/A$ has a right adjoint.
\end{thm}
As a consequence we get the following.

\begin{cor} If $p: A \to B$ is a covering morphism of
crossed complexes, then $p^* :
{\Crs}/B \to {\Crs}/A$ preserves all colimits.
\label{prescolimit}
\end{cor}

\section{Coverings of free crossed complexes} \label{free}
We recall here a definition from \cite{BH-81}. A {\em free basis}
for a  crossed complex $C$ consists of subgraphs $X_n$ of $C_n$
for all $ n \ge 1$  such that $C_1$ is the  free groupoid on
$X_1$, $C_2$ is the  free crossed $C_1$-module on the restriction
$\delta _2 ': X_2 \to C_1$, and for $n\geq3$, $C_n$ is the free
$\pi_1{C}$-module on $X_n$.

Following \cite{BH-91} we write ${\mathbb  C}(n)$ for the crossed
complex freely generated by one generator $c_n$ in dimension $n$.
So ${\mathbb  C}(0)$ is the singleton set $\{1\}$; ${\mathbb
C}(1)$ is the groupoid $\cal I$  which has two objects 0, 1 and
non-identity elements $c_1 : 0 \to 1$ and $c_{1}^{-1} : 1 \to 0$;
and for $n \geq 2, {\mathbb  C}(n)$ is  in dimensions $n$ and
$n-1$ an infinite cyclic group with generators $c_n$ and $\delta
c_n$ respectively, and is otherwise trivial. Thus if $C$ is a
crossed complex, then an element $c \in C_n$ is completely
specified by a morphism $\hat{c} : {\mathbb  C}(n) \to C$ such
that $\hat{c}(c_n)=c$, and $\delta(c)= \hat{c}(\delta c_n)$.

Let ${\mathbb  S}(n-1)$ be the subcomplex of ${\mathbb  C}(n)$
which agrees with ${\mathbb  C}(n)$ up to dimension $n-1$ and is
trivial otherwise. If ${\bf E}^n$ and ${\bf S}^{n-1}$ denote the
skeletal filtrations of the standard $n$-ball and $(n-1)$-sphere,
where $E^0 =\{ 0\} , S^{-1}= \emptyset, E^1 =I=\{ 0,1\} \cup e^1,
S^0 = \{0,1\} ,$ and for $n \ge 2, E^n = \{ 1\} \cup e^{n-1}\cup
e^n, S^{n-1} = \{ 1\} \cup e^{n-1}$, then it is clear that for all
$n \ge 0, {\mathbb  C}(n) \cong \pi{\bf E}^n$ and ${\mathbb
S}(n-1) \cong \pi{\bf S}^{n-1}.$

We now model for crossed complexes the process for spaces known as
attaching cells. Let $A$ be any crossed complex. A sequence of
morphisms $j_n : C^{n-1} \to C^n$ may be defined with $C^0 =A$ by
choosing any family of morphisms ${\mathbb  S}(m_{\lambda} -1) \to
C^{n-1}$ for any $\lambda \in \Lambda_n$ and any $m_{\lambda}$,
and forming the pushout
\begin{equation}{
\diagram
    \coprod_{\lambda \in  \Lambda_n}{\mathbb  S}(m_{\lambda}-1) \dto
\rto     & C^{n-1} \dto \\
    \coprod_{\lambda \in  \Lambda_n}{\mathbb  C}(m_{\lambda}) \rto
             & C^n.
\enddiagram } \label{freeconstruct} \end{equation} Let $C=
\mbox{colim}_{n}C^n$, and let $j:A \to C$ be the canonical
morphism. The morphism $j:A \to C$ is called a {\em relatively
free crossed complex morphism}. If $A$ is empty, then we call $C$
a {\em free} crossed complex.

The importance of the definition is as follows:
\begin{blank}\label{morfree}  If $C$ is a free crossed
complex on $X_*$, then a morphism $f:C \to D$ can be constructed
inductively provided one is given the values $f_n x\in D_n, x \in
X_n, n \ge 0$ provided the following geometric conditions are
satisfied: (i) $\delta^{\alpha} f_1x = f_0 \delta ^{\alpha }x, x
\in X_1, \alpha =0,1$; (ii) $ \beta f_n(x)= f_0(\beta x), x \in
X_n, n \ge 2$; (iii) $ \delta_n f_n(x) = f_{n-1}\delta_n (x), x
\in X_n, n \ge 2.$
\end{blank}

Notice that in (iii), $f_{n-1}$ has to be defined on
all of $C_{n-1}$ before this condition can be verified.

We now show that freeness can be lifted to covering crossed
complexes.

\begin{thm} \label{T:freecov}
Suppose given a pullback square of crossed complexes $$
\sqdiagram{\wt{A}} {\bar{j}}{\wt{C}}{p'}{p}{A}{j}{C}  $$ in which
$p$ is a covering morphism and $j : A \to C$ is relatively free.
Then $\bar{j} : \wt{A} \to \wt{C}$ is relatively free. \end{thm}

\begin{proof}  We suppose given the sequence of diagrams \ref{freeconstruct}.
Let $\hat{C}^n = p^{-1} (C^n). $ By corollary \ref{prescolimit},
 the following diagram is a pushout:
$$
\diagram
   p^* \left( \coprod_{\lambda \in  \Lambda_n}{\mathbb  S}(m_{\lambda}-1)\right) \dto \rto     & \hat{C}^{n-1} \dto \\
    p^*\left( \coprod_{\lambda \in  \Lambda_n}{\mathbb  C}(m_{\lambda}) \right) \rto               & \hat{C}^n.
\enddiagram
$$  Since $p$ is a covering morphism, we can write
 $  p^*\left( \coprod_{\lambda \in
\Lambda_n}{\mathbb  C}(m_{\lambda}) \right)  $ as
$  \coprod_{\lambda \in  \wt{\Lambda}_n}{\mathbb  C}(m_{\lambda})  $
for a suitable
 $\wt{\Lambda}_n$. This completes the proof.
\end{proof}
\begin{cor} \label{freecov}
Let $p : \wt C \to C$ be a covering morphism of crossed complexes.
If $C$ is free on $X \sast$, then $\wt{C}$ is  free on $p^{-1} (X
\sast)$.
\end{cor}

A similar result to Corollary \ref{freecov} applies in the
$m$-truncated case.

The significance of these results is as follows. We start with an
$m$-truncated free crossed resolution $C$ of a group $G$, so that
we are given $\phi: C_1 \to G$, and $C$ is free on $X \sast$,
where $X_n$ is defined only for $n \le m$. Our extension process
of section \ref{homotopies} will start by constructing the
universal cover $p: \wt{C}\to C$ of $C$; this is the covering
crossed complex corresponding to the universal covering groupoid
$p_0:\wt{G}\to G$. By the results above, $\wt{C}$ is the  free
crossed complex on $p\io(X\sast)$. It also follows from
Proposition \ref{isocov} that the induced morphism $\wt{\phi}:
\wt{C}\to \wt{G}$ makes $\wt{C}$ a free crossed resolution of the
contractible groupoid $\wt{G}$. Hence $\wt{C}$ is an acyclic and
hence, since it is free, also a contractible crossed complex.

We now see the general context for the diagram \eqref{maindiag} of
section \ref{compute} and the exposition there.

\section{Homotopies}
\label{homotopies}

We follow the conventions for homotopies in
\cite{Br-Hi2}. Thus a homotopy $  f ^0\simeq f $ of morphisms
$f^0,f : C \to D$ of crossed complexes is a pair $(h,f)$ where
$h$ is a family of functions $h_n : C_n \to D_{n+1} $ with the
following properties, in which $\beta c$ for $c \in C $ is $c,$
if $c \in C_0, $ is $\delta^1 c$, if $c \in C_1,$ and is $x$ if
$c \in C_n(x), \; n \ge 2.$ So we require \cite[(3.1)]{Br-Hi2}:
\begin{alignat}{2}
\beta h_n(c) &= \beta f(c) &\qquad &\mbox{for all } c \in C;\label{hom1} \\
h_1(cc') &= h_1(c)^{fc'}\;  h_1(c')&\qquad &
\mbox{if } c,c' \in C_1 \mbox{ and } cc' \mbox{ is defined;} \label{hom2}\\
h_2(cc') &= h_2(c)+ h_2(c') &\qquad &
\mbox{if } c,c' \in C_2 \mbox{ and } cc' \mbox{ is defined;} \label{hom22}\\
 h_n(c+c') &= h_n(c) + h_n(c') &\qquad & \mbox{if } c,c' \in
C_n,\; n \ge 3 \mbox{ and } c+c' \mbox{ is defined;}\label{hom3} \\
h_n(c^{ c_1} )&= (h_nc)^{fc_1} &\qquad & \mbox{if }  c \in C_n,
 n \ge 2, \;
c_1 \in C_1, \mbox{ and } c^{c_1} \mbox{ is defined.}
\label{hom4} \end{alignat}  Then $f^0,f $ are related by
\cite[(3.14)]{Br-Hi2}\begin{equation} f^0(c) = \left\{
\begin{array}{lll} \delta^0 h_0 c & \mbox{if} & c \in C_0, \\
(h_0 \delta^0 c)(fc)({\delta}_2h_1c)(h_0 \delta^1c)^{-1} & \mbox{if}
& c \in C_1, \\
\{( fc)( h_{1}{\delta}_2 c)( {\delta}_{3}h_2c) \}
^{(h_0\beta c)^{-1}} & \mbox{if} & c \in C_2, \\
 \{ fc + h_{n-1}{\delta}_n c + {\delta}_{n+1}h_nc \}
^{(h_0\beta c)^{-1}} & \mbox{if} & c \in C_n, \; n \ge 3.
\end{array}    \right. \label{homreln} \end{equation}

The following is important for our computations. We saw in
\ref{morfree} that a morphism is specified by its values on a
graded set of free generators. We now show that the same is true
for homotopies.

\begin{blank} \label{spechomot}
If $C$ is a free crossed complex on a generating family $X_n, n
\ge 0$, then a homotopy $(h,f):f^0 \simeq f: C \to D$ is specified
by the values $fx \in D_n, hx \in D_{n+1}, x \in X_n, n \ge 0$
provided only that the following geometric conditions hold:
\begin{equation}\begin{split}\delta^0fx&=f\delta^0x,\delta^1fx=f\delta^1x,
 x \in X_1,  \delta fx
=f \delta x, x \in X_n, n \ge 2,\\
\beta
fx&=f \beta x, x \in X_n, n \ge 1, \beta hx= \beta fx, x \in X_n, n
\ge 0. \end{split} \end{equation}
\end{blank}
\begin{proof}
All but the last condition are those given for the
construction of $f$ in \ref{morfree}. The final fact we need is
that for $n \ge 2$ the $f_1$-morphism $h_n$ is defined by its values on the
generators in $X_n$,  and this is standard.
\end{proof}

This result is another aspect of the facts that a homotopy $C \to
D$ can also be regarded as a morphism $\mathbb{C}(1) \otimes C \to
D$, where the tensor product is defined in \cite{Br-Hi2}, and the
tensor product of free complexes is free is proved in
\cite{BH-91}.

From this we can deduce formulae for a retraction. Suppose then
in the above formulae we take $C=D, \; f^0= 1_C, \; f=0$ where $0$ denotes
the constant morphism on $C$ mapping everything to a base point 0.
Then the homotopy $h : 1 \simeq 0$ must satisfy

\begin{alignat}{2}
\beta h_n c &= 0 \qquad & \text{if }c & \in C  , \label{retr1}\\
\delta^0 h_0 c & = c  \qquad & \text{if }c &  \in C_0 , \label{retr2}\\
{\delta}_2h_1c & = (h_0 \delta^0 c)\io c\, (h_0 \delta^1c)
\qquad & \text{if } c & \in C_1,  \label{retr3}\\
{\delta}_{3}h_2c & = (h_1{\delta}_2 c)\io  c^{h_0\beta c}
\qquad & \text{if }c & \in C_2, \label{retr32}\\
{\delta}_{n+1}h_nc & = -h_{n-1}{\delta}_n c+  c^{h_0\beta c}
\qquad & \text{if }c & \in C_n, n \ge 3 , \label{retr4}\\
h_n(c^{ c_1} )&= (h_nc) \qquad & \text{ if }  c &\in C_n,
 n \ge 2, \;
c_1 \in C_1, \mbox{ and } c^{c_1} \mbox{ is defined.}
\label{retr5}
 \end{alignat}
Further, in this case $h_1$ is a morphism by (\ref{hom2})
and for $n \ge 2$, $h_n$
 is by (\ref{hom3}) a  morphism which by (\ref{hom4})
trivialises the operations of $C_1$.  All these conditions are necessary
and sufficient for $h$ to be a contracting homotopy.

An {\em $m$-truncated crossed complex} $C$ is a crossed complex as
earlier except that $C_n$ and $\delta_n$ are defined only for $n
\le m$. Similarly, for a contracting homotopy $h$ of an
$m$-truncated crossed complex $C$, we have $h_n$ defined only for
$n < m$ and the above conditions hold where they make sense.

Our main result is now rather formal and  straightforward to
prove. It is to extend the pair $(C,h)$ of a  partial free crossed
resolution $C$ and partial contracting homotopy $h$ of the
universal cover of $C$ by one step. Hence the process can be
continued indefinitely.

\begin{thm} \label{extend}
Let $m \ge 1 $ and let   $C$ be an $m$-truncated free crossed
resolution of a group $G.$ Let $p:\wt{C} \to C$ be the universal
cover of $C$ so that $\wt{C}$ is an $m$-truncated free crossed
resolution of the universal covering groupoid $p_0: \wt{G}\to G$
of $G$. Let $h$ be a partial contracting homotopy of $\wt{C}.$
Suppose also that $\wt{C}_m$ is free on $\wt{X}_m.$

Let $\wt{X}_m \to X_{m+1}, \tilde{x} \mapsto x$, be a
 bijection to a set $X_{m+1}$ disjoint from $\wt{X}_m$.
Define an extension $e(C)$ of
$C$ to an $(m+1)$-truncated free crossed complex  as follows:\\
For $m=1$, let $C_2=e(C)_2$ be the free crossed $C_1$-module on
$X_2$ with $\delta_2 :X_2\to C_1  $ given by
\begin{alignat}{2}
{\delta}_2x  &= p_1\left((h_0 \wt{\delta}^0  \tilde{x})\io  \, \tilde{x}\,
\, (h_0 \wt{\delta}^1\tilde{x})\right), &\tilde{x}  &\in \wt{X}_1.  \label{ext1}
\\\intertext{For $m=2$ let $C_3=e(C)_3$ be the free $G$-module on $X_3$ with
$\delta_3 :e(C)_3\to C_2 $ defined by }
{\delta}_3 x  &= p_2 \left((h_1\wt{\delta}_2 \tilde{x})\io \; \tilde{x}
^{h_0\beta \tilde{x}}\right), & \tilde{x} &\in \wt{X}_2.\label{ext2} \\
\intertext{For $m\ge 3$ let $C_{m+1}=e(C)_{m+1}$ be the free $G$-module on $X_{m+1}$ with
$\delta_{m+1} :C_{m+1}\to C_{m}  $ defined by}
{\delta}_{m+1}x & =
p_m \left(-h_{m-1}\wt{\delta}_m \tilde{x}
+\tilde{x}^{h_0\beta \tilde{x}}\right),& \;
 \tilde{x} &\in \wt{X}_m.\label{extm}
\end{alignat}
Let $e(p): e(\wt{C})\to e(C)$ be the induced covering morphism,
extending $p$ by $p_{m+1}: \wtC_{m+1} \to C_{m+1}$.

Define $h_m:  \wt{C}_m \to \wt{C}_{m+1}$ on the basis $\wt{X}_m$ by
$h_m(\tilde{x})= (1,x). $ Then this extension $e(h)$ of $h$
is a contracting homotopy of $ e(\wt{C})$. Hence  $e(C)$
is an $(m+1)$-truncated free crossed resolution of $G.$

If further there is a subset $Y$ of $X_{m+1}$ such that
$\delta_{m+1}Y$ also generates $\ker \delta _{m}$, and a
retraction $\xi : C_{m+1}\to C_{m+1}(Y)$ is given such that
$\delta_{m+1}\xi (x)= \delta _{m+1}(x)$ for all $x \in X_{m+1}$,
and $\xi$ is a $G$-morphism for $m \ge 2$, and a crossed
$C_1$-morphism for $m=1$, then we may replace $C_{m+1}$ by
$C_{m+1}(Y)$ and $h_m$ by $\xi h_m$ to again get an extension of
the pair $(C,h)$ by one step.

\end{thm}
\begin{proof}
The fact that we have a contracting homotopy is immediate from the
definitions. It follows that $ e(\wt{C})$ is exact, and so $e(C)$
is aspherical with $\pi_1(e(C)) = G. $
\end{proof}
\begin{cor}\label{syzygies} Under the assumptions of Theorem {\em \ref{extend}},
 if $m=1$ then $\Ker \phi : C_1 \to G$ is generated as a normal subgroup of
$C_1$ by the elements:
 \begin{alignat}{3}
&p_1\left((h_0 \wt{\delta}^0  \tilde{x})\io  \, \tilde{x}\,
\, (h_0 \wt{\delta}^1\tilde{x})\right), &\;\tilde{x}  &\in \wt{X}_1.&&\\
\intertext{For $m \ge 2$,
$\Ker (\delta_m: C_m \to C_{m-1})$ is generated as a $G$-module by the elements:}
&p_2 \left((h_1\wt{\delta}_2 \tilde{x})\io \; \tilde{x}
^{h_0\beta \tilde{x}}\right),& \; \tilde{x} &\in \wt{X}_2,&
\text{ } & \text{  if } m=2,\\
&p_m \left(-h_{m-1}\wt{\delta}_m \tilde{x}
+\tilde{x}^{h_0\beta \tilde{x}}\right),& \;
 \tilde{x} &\in \wt{X}_m, &\text{ } & \text{ if }m \ge 3.
\end{alignat}
\end{cor}

We have now finally justified the process set out in section
\ref{compute} and illustrated with an example in section  \ref{syzygs3}.

In papers to follow we will give implementations of these methods,
and so  a wider range of calculations.

\section{Examples}
\label{moreres}
\subsection{The standard crossed resolution of a group}
\label{standres}
The standard crossed resolution of a group was defined by
Huebschmann in \cite{Hue} and applied also in for example
\cite{B-P,T}. Here we show how this resolution  arises from our
procedure.

We start with a group $G$ and let $C_1=F(G)$, the
free group on the set $G$, with generators written $[a], \; a \in
G$. Let $\phi: C_1 \to G$ be the canonical morphism. This has a
section $\sigma: G \to F(G), \, a \mapsto [a],a \ne 1,  1 \mapsto 1$.
This defines $h_0: G \to \wtC _1, \, a \mto (a,[a] \io)$.

The Cayley graph of this presentation has arrows $(a,[b]): a \to
ab$ so that $h_0(a)\io (a,[b])h_0(ab)= (1,[a][b][ab]\io)$. So we
may take $C_2$ to be the free crossed
$C_1$-module on elements $[a,b]$ and define $\delta_2:C_2 \to
C_1$ by $$\delta_2[a,b]= [a][b][ab]\io .$$
Then in the universal cover we can define $h_1: \wtC_1 \to
\wtC_2(1)$ by $(a,[b]) \mapsto (1,[a,b])$.

\begin{thm}
There is a free crossed $C_*(G)$ resolution of a group $G$ in
which $C_n(G)$ is free on the set $G^n$ with generators written
$[a_1,a_2,\ldots,a_n], a_i \in G$, with contracting homotopy of
the universal cover given by $(a,[a_1,a_2,\ldots,a_n])
 \mto (1,[a,a_1,a_2,\ldots,a_n])$, and boundary $\delta_n: C_n(G) \to C_{n-1}(G)$
 given by $\delta_2$ as above,
$$\delta_3[a,b,c]= [a,bc][ab,c]\io[a,b]\io [b,c]^{[a]\io},$$
and for $n \ge 4$
\begin{multline}
\delta_n[a_1, a_2, \ldots, a_n]=[a_2, \ldots,a_n]^{a_1\io} +
 \sum _{i=1}^{n-1}(-1)^{i}
[a_1,a_2,\ldots,a_{i-1}, a_ia_{i+1},a_{i+2}, \ldots, a_n] +\\
\qquad \qquad +(-1)^{n} [a_1,a_2, \ldots, a_{n-1}].
\end{multline}
\end{thm}
\begin{proof}
We first verify
\begin{align*}
\wtd_3h_2(a,[b,c]) &= h_1\wtd_2(a,[b,c])\io (a,[b,c])^{(a,[a]\io)}
\\ &= h_1(1,a,[b][c][bc]\io)\io(1,[b,c]^{[a]\io}) \\
&= h_1((a,[b])(ab,[c])(abc,[bc]\io))\io(1,[b,c]^{[a]\io})\\
&= h_1((a,[b])(ab,[c])(a,[bc])\io )\io(1,[b,c]^{[a]\io})\\
&= (1, [a,bc][ab,c]\io[a,b]\io[b,c]^{[a]\io}).
\end{align*}
 In order to have a contracting homotopy we require for $n \ge 3$
\begin{align*} &
\wtd_{n+1}h_n(a_1,[a_2,\ldots,a_{n+1}])\\ &=
-h_{n-1}\wtd_n(a_1,[a_2,\ldots,a_{n+1}]) +
(1,[a_2,\ldots,a_{n+1}]^{a_1\io})\\ &=
(1,[a_2,\ldots,a_{n+1}]^{a_1\io}+ \sum _{i=1}^{n}(-1)^{i}
[a_1,\ldots, a_ia_{i+1}, \ldots, a_{n+1}] +(-1)^{n+1} [a_1,a_2,
\ldots, a_n]).
\end{align*}
This completes the proof that the family $h_n$ give a contracting homotopy
and so that $C_*(G)$ is a resolution.

\end{proof}

\subsection{A small crossed resolution of finite cyclic groups}
This is the resolution given in \cite{BW}. Here we shall describe
its universal cover and a contracting homotopy. We would like to
thank A. Heyworth for discussions on this section.

We write $\cc$ for the (multiplicative) infinite cyclic group with
generator $x$, and $\cc_r$ for the finite cyclic group of order
$r$ with  generator $t$. Let $\phi: \cc \to \cc_r$ be the morphism
sending $x$ to $t$. We show how the inductive procedure given
earlier recovers the small free crossed resolution of $\cc_r$
together with a contracting homotopy of the universal cover.

Let $p_0: \wt{\cc}_r \to \cc_r$ be the universal covering
morphism, and let $p_1: \wt{\cc} \to \cc$ be the induced cover of
$\cc$. Then $\wt{\cc}$ is the free groupoid on the Cayley graph
$\wtX$ pictured as follows:

\begin{equation*}\xymatrixcolsep{3pc}
\xymatrix{
1^{\rule{0em}{1ex}} \ar@<-0.2ex> [r] _{(1,x)} & t^{\rule{0em}{1ex}}
 \ar@<-0.2ex> [r] _{(t,x)}&
 t^2 \ar@<-0.2ex> [r] _{(t^2,x)} &
\ar @{.}@<-0.2ex> [r] & \ar@<-0.2ex> [r] & t^{n-2} \ar @<-0.2ex>[r] _{(t^{n-2},x)}
& t^{n-1} \ar @/_2pc/
[0,-6] _{(t^{n-1},x)}
}
\end{equation*}
A  section $\sigma: \cc_r \to \cc$ of $\phi$ is given by $t^i \mto
x^i, i=0,\ldots, r-1$, and this defines
$$ h_0: \cc_r
\to \wtF_1 , \; t^i \mto (t^i,x^{-i}).$$
It follows that for $ i=0,\ldots, r-1$
$$ h_0(t^i)\io (t^i,x)h_0(t^{i+1}) =
\begin{cases} (1,1) & \text{if } i \ne r-1, \\
              (1,x^r) & \text{if } i = r-1.
\end{cases} $$ So we take a new generator $x_2$ for $F_2$ with
$\delta_2 x_2=x_r$ and set
$$ h_1(t^i,x) =\begin{cases} (1,1) & \text{if } i \ne r-1, \\
              (1,x_2) & \text{if } i = r-1.
\end{cases} $$
Then for all $i =0, \ldots, r-1$ we have
$$ \wtd_2 h_1(t^i,x)= h_0(t^i)\io (t^i,x)h_0(t^{i+1}),$$
and it  follows that
\begin{align} h_1(t^i,x^r) &=h_1((t^i,x)(t^{i+1},x)\ldots
                    (t^{i+r-1},x))   \notag\\
            &= (1,x_2). \label{rhomot}
\end{align}
Hence \begin{align*}
-h_1 \wtd_2 (t^i,x_2) + (t^i,x_2).x^{-i}&=
(1,-x_2)+(1,x_2.t^{-i})\\
&= (1,x_2.(t^{r-i}-1))
\end{align*}
This gives us 0 for $i=0$, and $(1,x_2.(t-1))$ for $i=r-1$. Let
$N(i)= 1+t+ \cdots+t^{i-1}$, so that $t^{r-i}-1=(t-1)N(r-i)$ for
$i=1,\ldots,r-1.$ Hence we can take a new generator $x_3$ for
$F_3$ with $\delta_3x_3= x_2.(t-1)$ and define
$$h_2(t^i,x_2) =\begin{cases} (1,0) & \text{if } i =0, \\
              (1,x_3.N(r-i)) & \text{if } 0< i \le r-1.
\end{cases} $$
Now we find that if we evaluate
$$-h_2\wtd_2(t^i,x_3) + (1,x_3.t^{-i})= -h_2((t^{i-1},x_2).t+(t^i,x_2))+  (1,x_3.t^{-i})$$
we obtain for $i=0$
$$ -h_2(t^{r-1},x_2) +(1,x_3) =(1,0),$$
for $i=1$
$$0+ h_2(t,x_2) +(1,x_3.t^{r-1}) =
(1,x_3.(N(r-1)+t^{r-1}))=(1,x_3.N(r))$$
and otherwise
$$(1,x_3(-N(r-i+1)+N(r-i)+t^{r-i}) )=(1,0).$$
Thus we take a new generator $x_4$ for $F_4$ with $\delta_4x_4=
x_3.N(r)$ and $$h_3(t^i,x_3)=\begin{cases} (1,x_4) & \text{if } i =1, \\
              (1,0) & \text{otherwise.}
\end{cases} $$
Then \begin{align*}
-h_3\wtd_4(t^i,x_4)+(1,x_4.t^{-i}) &= -h_3(t^i,x_3.N(r))
+(1,x_4.t^{-i})\\
&=-h_3(1,x_3.N(r).t^{-i})+(1,x_4.t^{-i})\\
&= (1,x_4. (t^{r-i}-1)).
\end{align*}
Thus we are now in a periodic situation and we have the theorem:
\begin{thm}
A free crossed resolution $F_*$ of $\cc_r$ may be taken to have
single free generators $x_n$ in dimension $n \ge 1$ with $\phi
(x_1)= t$, and $$
 \delta_n(x_n) = \begin{cases}
        x_1^r &\text{if } n=2, \\
        x_{n-1}.(t-1)&\text{if } n>1, n \text{ odd}, \\
 x_{n-1}.N(r)&\text{if } n>2, n \text{ even}.
 \end{cases} $$
A contracting homotopy on $\wtF_*$ is given by $h_0(t^i)=
(t^i,x_1^{-i}) $ and for $n>1$
$$h_n(t^i,x_n) = \begin{cases}
(1,0)&\text{if } n=1, i \ne r-1, \\
(1,x_2)&\text{if } n=1, i = r-1, \\
(1,0)&\text{if } n \text{ even, }n\ge 2, i =0, \\
(1,x_{n+1}.N(r-i))&\text{if } n \text{ even, }n\ge 2, 0<i \le r-1,\\
(1,0)&\text{if } n \text{ odd, }n>1, i \ne 1, \\
(1,x_{n+1})&\text{if } n \text{ odd, }n>1,  i=1.
\end{cases}$$
\end{thm}


\begin{thebibliography}{99}
\bibitem{Ba-Pr}{\sc Baik, Y.G., and Pride, S.J.}, `Generators of the
second homotopy module of presentations arising from group
constructions', Preprint, Glasgow 1992.
\bibitem{Blakers}{\sc Blakers, A.L.}, `Some relations between
homology and homotopy groups', {\em Annals of Math.}, 49 (1948) 428-461.
\bibitem{B-70}{\sc  Brown, R.},  `Fibrations of groupoids', {\em J.
Algebra} 15 (1970) 103-132.
\bibitem{B-88}{\sc  Brown, R.},  {\em Topology: A geometric account
of general topology,
homotopy types and the fundamental groupoid} (Ellis Horwood, 1988).
\bibitem{B-Irish}{\sc  Brown, R.},
`Higher order symmetry of graphs', {\em Bull. Irish Math.
Soc.} 32 (1994) 46-59.
\bibitem{B-94} {\sc  Brown, R.}, `Homotopy theory, and change of base
for groupoids and multiple groupoids', {\em Applied categorical structures},
4 (1996) 175-193.
\bibitem{B-gre}{\sc Brown, R.}, `Groupoids and crossed objects
in algebraic topology', {\em Homotopy,  homology and applications},
(to appear).
\bibitem{BHK} {\sc Brown, R., Heath, P.R.  and
Kamps, H.}, `Coverings of groupoids
and  Mayer-Vietoris type sequences', {\em Categorical Topology,
Proc. Conf.  Toledo, Ohio, 1983}.  (Heldermann Verlag, Berlin)
(1984) 147-162.
\bibitem{BH1} {\sc Brown, R. and Higgins, P.J.}, `On the
connection between the second relative homotopy groups of some
related spaces',  {\em Proc. London Math. Soc.}, (3) 36 (1978)
193-212.
\bibitem{BH2} {\sc Brown, R. and Higgins, P.J.}, `On the algebra
of cubes', {\em J. Pure Appl.  Algebra}, 21 (1981) 233-260.
\bibitem{BH3} {\sc Brown, R. and Higgins, P.J.}, `Colimit theorems for
 relative homotopy  groups', {\em J. Pure Appl. Algebra}, 22 (1981) 11-41.
\bibitem{BH-81} {\sc Brown, R.  and  Higgins, P.J.}, `Crossed complexes and
 non-abelian extensions',  In {\em Proc. International Conference on Category Theory:
 Gummersbach, 1981}, Lecture Notes in Math. vol. 962 (Springer-Verlag, 1982), 39-50.
\bibitem{Br-Hi2} {\sc Brown, R.  and  Higgins, P.J.}, `Tensor
products and homotopies for  $\omega$-groupoids and crossed
complexes',  {\em  J. Pure Appl. Algebra}, 47 (1987), 1-33.
\bibitem{BH-90}{\sc Brown, R. and Higgins, P.J.},  `Crossed complexes and chain
complexes with  operators', {\em Math. Proc. Camb. Phil. Soc.}, 107 (1990) 33-57.
\bibitem{BH-91} {\sc Brown, R. and Higgins, P.J.}, `The classifying space
 of a crossed complex',  {\em Math. Proc. Camb. Phil. Soc.}, 110 (1991), 95-120.
\bibitem{B-Hu} {\sc Brown, R. and Huebschmann, J.},  `Identities
among relations', in  {\em Low dimensional topology},  Ed.
R.Brown and T. L. Thickstun,  London Math. Soc. Lecture Notes
46,  Cambridge University Press (1982) 153-202.
\bibitem{B-P}{\sc Brown, R. and Porter, T.}, `On the Schreier theory of
non-abelian extensions: generalisations and computations',
{\em Proceedings Royal Irish Academy}, 96A (1996) 213-227.
\bibitem{BW}  {\sc Brown, R. and Wensley, C.D.},    `On finite
induced crossed modules, and the homotopy 2-type  of mapping
cones', {\em Theory and Applications of Categories}, 1 (1995)
54-71.
\bibitem{BW2}  {\sc Brown, R. and Wensley, C.D.},   `Computing crossed
modules induced by an
inclusion of a normal subgroup, with applications to homotopy
2-types', {\em Theory and Applications of Categories} 2 (1996) 3-16.
\bibitem{Fox} {\sc Fox, R.H.},  `Free differential calculus I. Derivations
in the free group ring', {\em Annals of Math.} 57 (1953) 547-560.
\bibitem{Gro} {\sc Groves, J.R.J.}, `An algorithm for computing
homology groups', {\em J. Pure Appl. Algebra} 194 (1997) 331-361.
\bibitem{Gru1}{\sc Gruenberg, K.W.}, `Resolution by relations',
{\em J. London Math. Soc.}, 35 (1960) 481-494.
\bibitem{He} {\sc Heyworth, Anne}, {\em  Applications of Rewriting Systems and
   Gr\"obner Bases to Computing Kan Extensions  and Identities Among
    Relations}, University of Wales Ph.D. thesis (1998).
\\ {\tt http://xxx.soton.ac.uk/abs/math.CT/9812097}
\bibitem{HeR} {\sc Heyworth, Anne and Reinert, Birgit},
`Applications of Gr\"{o}bner bases to group rings', (in
preparation).
\bibitem{HeW} {\sc Heyworth, Anne and Wensley, C.}, `Logged
Knuth-Bendix procedures and identities among relations', (in
preparation).
\bibitem{Higg-Pres}{\sc Higgins, P.J.}, `Presentations of
 groupoids, with applications to  groups', {\em Proc. Camb.
 Phil. Soc.},  60 (1964) 7-20.
\bibitem{Hi} {\sc Higgins, P.J.},  {\em  Categories and
groupoids} Van Nostrand, (1971).
\bibitem{HMS}    {\sc Hog-Angeloni, C., Metzler, W. and Sieradski,
A.J. (Editors)},    {\em Two-dimensional homotopy and
combinatorial group theory},    London Math. Soc. Lecture Note
Series 197,    Cambridge University Press, Cambridge (1993).
\bibitem{Howie} {\sc Howie, J.}, `Pullback functors and  crossed complexes',
{\em Cahiers Topologie G\'{e}om. Diff\'{e}rentielle Cat\'{e}goriques}
20 (1979) 281-295.
\bibitem{Hu} {\sc  Hu, S.T.}, {\em Homotopy theory} (Academic Press, 1959).
\bibitem{Hue} {\sc Huebschmann}, `Crossed $n$-fold extensions and
cohoology', {\em Comm. Math. Helv.}, 55 (1980) 302-314.
\bibitem{JLS} {\sc Johansson, L., Lambe, L., and Sk\"oldberg, E.}, `Normal forms
and iterative methods for constructing resolutions', Preprint Stockholm,
1998.
\bibitem{P} {\sc Pride, S.J.},  `Identities among relations', in
{\em Proc. Workshop on Group Theory from  a  Geometrical
Viewpoint},  eds. E.Ghys, A.Haefliger, A. Verjodsky,
International Centre  of  Theoretical  Physics, Trieste, 1990,
World Scientific (1991) 687-716.
\bibitem{GAP}  {\sc Sch{\accent127 o}nert, M. et al},  GAP :
Groups, Algorithms, and Programming,   Lehrstuhl D f{\accent127
u}r Mathematik, Rheinisch  Westf{\accent127 a}lische Technische
Hoch\-schule, Aachen, Germany, fourth edition, 1997.
\bibitem{T} Tonks, A.P., {\em  Theory and applications of
 crossed complexes}, PhD thesis, University of Wales, Bangor
 (1993)\\ { \tt  http://www.bangor.ac.uk/ma/research/tonks}.
\bibitem{W-A} {\sc Wensley, C.D. and Alp, M.}, {\em XMOD}, Share package
distributed by the GAP Council, 1997, Ch 73 of \cite{GAP}.
\bibitem{W1} {\sc  Whitehead, J.H.C.},  `Combinatorial homotopy
II',  {\em  Bull. Amer. Math. Soc.} 55 (1949) 453-496.

\end{thebibliography}
\end{document}